\documentclass[11pt]{article}
\usepackage[english]{babel}
\usepackage[T1]{fontenc}
\usepackage[utf8]{inputenc}
\usepackage{amsmath,amssymb,dsfont, float, bm}
\usepackage{amsmath,amssymb,amsfonts,amsthm,mathtools,mathrsfs,booktabs}
\usepackage{tikz,float}
\usepackage{pgfplots}
\pgfplotsset{compat=1.16}
\usepackage{authblk}

\usepackage{natbib}

\usepackage[colorlinks=true, citecolor=blue, urlcolor=blue]{hyperref}
\usepackage{newtxtext,newtxmath} % Times Roman fonts (text and math)
\usepackage{mhchem} % for \ce macro
\usepackage{multirow}
\usepackage{bookmark}
\usepackage{booktabs}

\newcommand{\E}{\mathbb{E}}

\newcommand{\Prob}{\mathbb{P}}
\newcommand{\R}{\mathbb{R}}

\newcommand{\argmin}{\operatorname*{arg\,min}}

\hypersetup{                 
colorlinks=true,
linkcolor=blue,
citecolor=blue,
anchorcolor=blue,
urlcolor=blue,
linktoc=page,
}
\usepackage[a4paper,margin=2.5cm]{geometry}

\newtheorem{theorem}{Theorem}%[section]

\newtheorem{example}{Example}[section]

\newtheorem{lemma}{Lemma}[section]
\newtheorem{proposition}{Proposition}[section]
\newtheorem{remark}{Remark}

\begin{document}

\title{\bf 
Jointly optimized fractional pairwise log-difference estimators for the two-parameter Weibull distribution
}

%\author{}
 \author[]{Roberto Vila}%\thanks{Corresponding author:  rovig161@gmail.com}}
 \author[]{\, Felipe Quintino \thanks{Corresponding author: felipes.quintino2@gmail.com}}
 %\author[]{Helton Saulo }
 \affil[]{Department of Statistics, University of
 	Bras\'ilia, 70910-900, Bras\'ilia, Brazil}
 \date{}                     %% if you don't need date to appear
\setcounter{Maxaffil}{0}
\renewcommand\Affilfont{\itshape\small}

%\date{}
\maketitle
\begin{abstract}
In this paper, we develop explicit estimators for the shape and scale parameters of the two-parameter Weibull distribution. The construction relies on fractional absolute moments of pairwise logarithmic differences and a weighted combination of pairwise logarithmic minima and maxima. Exploiting the logistic distribution of the difference between two independent logarithms of Weibull random variables, we obtain a family of shape estimators indexed by a positive fractional order and corresponding scale estimators controlled by a weighting parameter.
We establish strong consistency and derive the joint asymptotic distribution of the proposed estimators. We also formulate a determinant-based joint optimization criterion and show that the optimal weight has an explicit expression for each fixed fractional order, reducing the problem to one-dimensional optimization. We discuss the connections with existing estimators based on U-statistics, order statistics, and fractional moments.
Finally, we present a numerical evaluation of the fractional pairwise log-difference (FPLD) estimator and compare it with existing estimators from the literature under both uncontaminated and contaminated scenarios. The numerical results suggest that the FPLD estimator exhibits robustness to data contamination.  
\end{abstract}

%{
%\hypersetup{linkcolor=black}
%\tableofcontents
%}

%------------------------------------------------
\section{Introduction}
%------------------------------------------------

The Weibull distribution is widely used in reliability
analysis, survival studies, engineering, and lifetime
modeling. Its flexibility is mainly due to the shape
parameter, which determines the behavior of the
hazard function, while the scale parameter controls
the magnitude of the observations.
The Weibull model has applications in several fields, such as in medicine \citep{matsushita1992lifetime}, risk insurance \citep{Gebizlioglu2011}, hydrology \citep{Clarke2002},  climatology \citep{Papalexiou2018}, finance \citep{mittnik1993modeling}, and engineering \citep{kundu2009estimation}. A complete characterization of the Weibull model and its applications can be found in the books of \cite{Johnson1994} and \cite{Murthy2004}.

Despite the simplicity of its density, the estimation
of its parameters is not always straightforward.
Maximum likelihood estimation is a standard
approach, but the likelihood equation for the shape
parameter does not generally admit a closed-form
solution. Consequently, numerical procedures are
usually required.

Several alternatives have been investigated. 
These include percentiles methods\citep{Dubey1967},  L-moments \citep{Hosking1990}, modified maximum likelihood \citep{Jokiel2024},
logarithmic moments, probability-weighted moments,
least-squares procedures, moment-based,
and estimators based on order statistics. A comparative study of several classical methods was presented by
\citet{Teimouri2013}.
Estimators based on U-statistics were subsequently investigated by \citet{Sadani2019}. Both the Sadani estimator and the earlier Menon estimator \citep{Menon1963} arise as special cases of the proposed FPLD estimator. 

More recent contributions have continued to
investigate explicit and computationally convenient
estimators. For example, \citet{Jokiel2024}
proposed estimators based on modifications of
maximum likelihood equations and nonparametric
Gini estimators. \citet{Karahasan2024} developed
closed-form estimators that can also be adapted to
doubly Type II censored samples. Earlier work on
Weibull estimation based on selected order
statistics includes \citet{Murthy1975}.

Fractional moments have also received renewed
attention. In particular, \citet{AlnssyanDar2026}
studied a maximum-entropy characterization of
the Weibull distribution through fractional moments
of the original random variable. This approach
differs from the one considered here, which is based
on fractional moments of pairwise logarithmic
differences.

The purpose of this paper is to investigate a class
of explicit estimators based on the distributional
properties of logarithmic differences. The central
observation is that, if $X_1$ and $X_2$ are
independent Weibull random variables with common
shape and scale parameters, then
$
\log(X_1)-\log(X_2)
$
has a logistic distribution whose scale depends
only on the shape parameter.

This result motivates a family of estimators indexed
by a positive fractional order. To estimate the scale
parameter, we consider a weighted combination
of pairwise logarithmic minima and maxima.
The two tuning parameters are subsequently
considered jointly through the asymptotic covariance
matrix of the resulting estimators.

The contributions developed in this paper are
threefold. First, we introduce a fractional-order
extension of a pairwise logarithmic-difference
estimator for the Weibull shape parameter.
Second, we construct a corresponding family of
scale estimators and establish their joint
asymptotic properties. Third, we derive an explicit
reduction of the joint determinant-based optimization
problem. Finally, an extensive numerical evaluation of the FPLD estimator is conducted. The study examines the impact of the fractional order, proposes a reference table for selecting the tuning parameters based only on the sample size, compares this approach with optimization-based tuning procedures, and evaluates the performance of the FPLD estimator against both closed-form and iterative competitors.

The paper is organized as follows. Section \ref{Model and proposed estimators}
introduces the model and the proposed estimators.
Section \ref{Consistency} establishes their strong consistency.
Section \ref{Joint asymptotic distribution} derives the joint asymptotic distribution.
Section \ref{Joint asymptotic optimization} develops the asymptotic optimization
criterion. 
Section \ref{Numerical investigation} presents the numerical investigations, and Section \ref{Concluding remarks} provides concluding remarks.

\section{Model and proposed estimators}
\label{Model and proposed estimators}

Let $X_1,\ldots,X_n$ be an independent and
identically distributed (i.i.d) sample from the two-parameter
Weibull distribution with density
\begin{equation}
	f(x;k,\lambda)
	=
	\frac{k}{\lambda}
	\left(\frac{x}{\lambda}\right)^{k-1}
	\exp\left\{-\left(\frac{x}{\lambda}\right)^k\right\},
	\quad x>0,
	\label{eq:weibull}
\end{equation}
where $k>0$ and $\lambda>0$ denote the shape and
scale parameters, respectively.

Define
\begin{equation}
	Y_i=k\log\left(\frac{X_i}{\lambda}\right),
	\quad i=1,\ldots,n.
	\label{eq:Gtransform}
\end{equation}
Since
\[
\left(\frac{X_i}{\lambda}\right)^k
\sim \operatorname{exp}(1),
\]
the random variable $Y_i$ has distribution function
\begin{equation}
	F_Y(y)=1-\exp(-{\rm e}^y),
	\quad y\in\mathbb{R},
	\label{eq:gumbelcdf}
\end{equation}
and density
\begin{equation}
	f_Y(y)={\rm e}^y\exp(-{\rm e}^y),
	\quad y\in\mathbb{R}.
	\label{eq:gumbelpdf}
\end{equation}
Thus, $Y_i$ follows the standard Gumbel
distribution of the minimum type, denoted by $Y_i\sim{\rm Gumbel}(0,1)$.

An important consequence of this representation
is the following result.

\begin{lemma}\label{1-lemma}
	Let $Y_1$ and $Y_2$ be independent random variables
	with distribution function given by
	\eqref{eq:gumbelcdf}. Then
	\[
	Y_1-Y_2\sim\operatorname{Logistic}(0,1).
	\]
	Consequently,
	\[
	\log(X_1)-\log(X_2)
	\sim\operatorname{Logistic}(0,1/k).
	\]
\end{lemma}
\begin{proof}
	Let $Z_1$ and $Z_2$ be independent standard
	exponential random variables. By
	\eqref{eq:Gtransform},
	$
	Y_1-Y_2=\log(Z_1/Z_2).
	$
	For $z\in\mathbb{R}$,
	\[
	\begin{aligned}
		\Prob(Y_1-Y_2\leqslant z)
		=\Prob(Z_1\leqslant {\rm e}^z Z_2)
		=\int_0^\infty
		(1-{\rm e}^{-{\rm e}^z t}){\rm e}^{-t}\,{\rm d}t
		=\frac{{\rm e}^z}{1+{\rm e}^z}.
	\end{aligned}
	\]
	This is the standard logistic distribution
	function. The second assertion follows from
	$
	\log(X_1)-\log(X_2)=(Y_1-Y_2)/k.
	$
\end{proof}

%------------------------------------------------
\subsection{Fractional pairwise estimator of the
	shape parameter}
%------------------------------------------------

%For $r>0$, define
%\begin{equation}
%	M_r=2\Gamma(r+1)\etaD(r),
%	\label{eq:Mr}
%\end{equation}
%where $\etaD$ denotes the Dirichlet eta function,
%\[
%\etaD(r)=
%\sum_{j=1}^{\infty}
%\frac{(-1)^{j-1}}{j^r},
%\quad r>0.
%\]
%Let 
%\[
%\etaD(r)=
%\sum_{j=1}^{\infty}
%\frac{(-1)^{j-1}}{j^r},
%\quad r>0,
%\]
%be the Dirichlet eta function.
The absolute moments of $L\sim\operatorname{Logistic}(0,1)$ satisfy
\begin{equation}
	\E[|L|^r]=2\Gamma(r+1)
	\sum_{j=1}^{\infty}
	\frac{(-1)^{j-1}}{j^r},
	\quad r>0.
	\label{eq:logisticmoment}
\end{equation}

Define the pairwise statistic
\begin{equation}
	U_{n,r}
	\coloneqq
	{1\over\binom{n}{2}}
	\sum_{1\leqslant i<j\leqslant n}
	|\log(X_i)-\log(X_j)|^r.
	\label{eq:Unr}
\end{equation}

Since  $X_1,\ldots,X_n$ are i.i.d random variables following
\eqref{eq:weibull} and $k[\log(X_1)-\log(X_2)]
\sim\operatorname{Logistic}(0,1)$  (Lemma \ref{1-lemma}),
\begin{equation}
	\E[U_{n,r}]
	=
	\E[|\log(X_1)-\log(X_2)|^r]
	\stackrel{\eqref{eq:logisticmoment}}{=}
	\frac{2\Gamma(r+1)\sum_{j=1}^{\infty}
		\frac{(-1)^{j-1}}{j^r}}{k^r}.
	\label{eq:Unrexpectation}
\end{equation}

This identity motivates the estimator
\begin{equation}
		\widehat{k}_r
		=
		\left[\frac{2\Gamma(r+1)\sum_{j=1}^{\infty}
			\frac{(-1)^{j-1}}{j^r}}{U_{n,r}}\right]^{1/r}.
	\label{eq:kestimator}
\end{equation}
The estimator is explicit and does not depend on
the unknown scale parameter.

\begin{example}
    For $r=1$, since $\sum_{j=1}^{\infty}
{(-1)^{j-1}}/{j}=\log(2)$ and $\Gamma(2)=1$, expression
\eqref{eq:kestimator} becomes
\begin{equation}
	\widehat{k}_1=
	\frac{2\log(2)}{
	{1\over\binom{n}{2}}
		\sum_{{1\leqslant i<j\leqslant n}}|\log(X_i)-\log(X_j)|
	}.
	\label{eq:kestimatorr1}
\end{equation}
\end{example}

\begin{remark}
    \begin{enumerate}
        \item This special case coincides with the U-statistic
estimator of the shape parameter proposed by
\citet{Sadani2019}.
For $r=2$, since
$$\sum_{1\leqslant i<j\leqslant n}\{\log(X_i)-\log(X_j)\}^2
=n\sum_{i=1}^n\{\log(X_i)-\overline{\log(X)}\}^2,$$
where
$$\overline{\log(X)}=n^{-1}\sum_{i=1}^n\log(X_i),$$ we have
$U_{n,2}=2S^2,$ with $S^2$ the sample variance of
$\log(X_1),\ldots,\log(X_n)$. 

\item Due to
$$2\Gamma(3)\sum_{j=1}^{\infty}(-1)^{j-1}/j^2=\pi^2/3,$$
\eqref{eq:kestimator} becomes
$\widehat{k}_2=\pi/(\sqrt{6}\,S)$, which is the
logarithmic-moment estimator of \citet{Menon1963}.
    \end{enumerate}
\end{remark}

%------------------------------------------------
\subsection{Pairwise estimator of the scale parameter}
%------------------------------------------------

For $a\in\mathbb{R}$, define
$
	h_a(x,y)=
	a
	\log(\min\{x,y\})
	+
	(1-a)
	\log(\max\{x,y\}).
$
The corresponding U-statistic is
\begin{equation}
	V_{n,a}
	\coloneqq 
	{1\over \binom{n}{2}}
	\sum_{1\leqslant i<j\leqslant n} h_a(X_i,X_j).
	\label{eq:Vnr}
\end{equation}

Let $\gamma$ denote the Euler--Mascheroni
constant. It is well known that, for standard minimum-type Gumbel random variables
$Y$, $Y_1$, and $Y_2$, with $Y_1$ and $Y_2$ independent, the following
properties hold:
\begin{align}\label{exp-Y}
\E[Y]=-\gamma,
\quad 
\E[\min\{Y_1,Y_2\}]=-\log(2)-\gamma,
\quad 
\E[\max\{Y_1,Y_2\}]=\log(2)-\gamma.
\end{align}

Therefore, since $X_1,\ldots,X_n$ are i.i.d random variables following
\eqref{eq:weibull} and $\log(X_i)=\log(\lambda)+Y_i/k$ with $Y_i\sim{\rm Gumbel}(0,1)$,
\begin{align}
	\E[V_{n,a}]
	=
	\E[h_a(X_1,X_2)]
	&=
	a
	\left\{
	\log(\lambda)
	+
	{\mathbb{E}[\min\{Y_1,Y_2\}]\over k}
	\right\}
	+
	(1-a)
	\left\{
\log(\lambda)
+
{\mathbb{E}[\max\{Y_1,Y_2\}]\over k}
\right\}
\nonumber
\\[0,2cm]
	&=
	\log(\lambda)+\frac{(1-2a)\log(2)-\gamma}{k}.
	\label{eq:Cexpectation}
\end{align}
%where
%\begin{equation}
%	\theta_a=-\gamma_E+(1-2a)\log(2).
%	\label{eq:theta}
%\end{equation}

Combining \eqref{eq:Cexpectation} with
\eqref{eq:kestimator}, we define
\begin{equation}
		\widehat{\lambda}_{a,r}
		=
		\exp\left\{
		V_{n,a}-\frac{(1-2a)\log(2)-\gamma}{\widehat{k}_r}
		\right\},
	\label{eq:lambdaestimator}
\end{equation}
where $\widehat{k}_r$ is as given in \eqref{eq:kestimator}.

Both estimators $\widehat{k}_r$ and $\widehat{\lambda}_{a,r}$ are available in explicit form. Since 
$$a\min\{u,v\}+(1-a)\max\{u,v\}=(u+v)/2+(1/2-a)|u-v|
~~\text{for}~ u,v\in\mathbb{R},$$
the statistic \eqref{eq:Vnr} satisfies
the identity
\[
V_{n,a}=\overline{\log(X)}+\left(\frac{1}{2}-a\right)U_{n,1},
\]
where $U_{n,1}$ is given by \eqref{eq:Unr} with $r=1$. Hence,
\[
\log(\widehat{\lambda}_{a,r})
=\overline{\log(X)}+\frac{\gamma}{\widehat{k}_r}
+\left(\frac{1}{2}-a\right)
\left\{U_{n,1}-\frac{2\log(2)}{\widehat{k}_r}\right\}.
\]
Two particular cases are described in the following remark. 

\begin{remark}
    \begin{enumerate}
        \item For $a=1/2$,
$$\widehat{\lambda}_{1/2,r}=\exp\{\overline{\log (X)}+\gamma/\widehat{k}_r\}~~r>0.$$
In particular,
$(\widehat{k}_2,\widehat{\lambda}_{1/2,2})$ is the
logarithmic-moment estimator of \citet{Menon1963}.
    \item By \eqref{eq:kestimatorr1}, $2\log(2)/\widehat{k}_1=U_{n,1}$,
so the last term vanishes when $r=1$. Thus
$\widehat{\lambda}_{a,1}=\exp\{\overline{\log( X)}+\gamma/\widehat{k}_1\}$
does not depend on $a$, and $(\widehat{k}_1,\widehat{\lambda}_{a,1})$
coincides with the estimator of \citet{Sadani2019} for every
$a\in\mathbb{R}$.
    \end{enumerate}
\end{remark}

\section{Consistency}
\label{Consistency}

We first establish the consistency of the proposed
estimators \eqref{eq:kestimator} and \eqref{eq:lambdaestimator}.

\begin{theorem}
	Let $X_1,X_2,\ldots$ be independent and identically
	distributed random variables following
	\eqref{eq:weibull}, with $k>0$ and $\lambda>0$.
	For any fixed $r>0$ and $a\in\mathbb{R}$,
	\[
	\widehat{k}_r\xrightarrow{\mathrm{a.s.}}k,
	\quad 
	\widehat{\lambda}_{a,r}
	\xrightarrow{\mathrm{a.s.}}\lambda,
	\]
	where $\xrightarrow{\mathrm{a.s.}}$ denotes almost sure convergence.
\end{theorem}
\begin{proof}
	The statistic $U_{n,r}$ in \eqref{eq:Unr} is a U-statistic of order
	two with kernel
	$
	|\log(x)-\log(y)|^r.
	$
	By \eqref{eq:Unrexpectation},
	\[
	\E[|\log(X_1)-\log(X_2)|^r]={2\Gamma(r+1)\sum_{j=1}^{\infty}
		\frac{(-1)^{j-1}}{j^r}\over k^r}<\infty.
	\]
	The strong law of large numbers for U-statistics \citep{Lee1990,Henze2024}
	therefore yields
	\[
	U_{n,r}\xrightarrow{\mathrm{a.s.}}\frac{2\Gamma(r+1)\sum_{j=1}^{\infty}
		\frac{(-1)^{j-1}}{j^r}}{k^r}.
	\]
	Since $2\Gamma(r+1)\sum_{j=1}^{\infty}
	{(-1)^{j-1}}/{j^r}>0$, the continuous mapping theorem
	gives
	\[
	\widehat{k}_r\xrightarrow{\mathrm{a.s.}}k.
	\]
	
	Similarly, $V_{n,a}$ in \eqref{eq:Vnr} is a U-statistic with an
	integrable kernel. Hence,
	\[
	V_{n,a}
	\xrightarrow{\mathrm{a.s.}}
	\E[V_{n,a}]
	\stackrel{\eqref{eq:Cexpectation}}{=}
	\log(\lambda)+\frac{(1-2a)\log(2)-\gamma}{k}.
	\]
	It follows that
	\[
	\begin{aligned}
		\log(\widehat{\lambda}_{a,r})
		=
		V_{n,a}-\frac{(1-2a)\log(2)-\gamma}{\widehat{k}_r}
		\xrightarrow{\mathrm{a.s.}}
%		\log(\lambda)+\frac{(1-2a)\log(2)-\gamma}{k}
%		-\frac{(1-2a)\log(2)-\gamma}{k}
%		=
		\log(\lambda).
	\end{aligned}
	\]
	Another application of the continuous mapping
	theorem completes the proof.
\end{proof}

\section{Joint asymptotic distribution}
\label{Joint asymptotic distribution}

We now derive the joint asymptotic distribution
of the proposed estimators \eqref{eq:kestimator} and \eqref{eq:lambdaestimator}.

For $y\in\mathbb{R}$ and $Y\sim{\rm Gumbel}(0,1)$, define
\begin{equation}
	v_r(y)
	\coloneqq
	\int_{-\infty}^{\infty}
	|y-t|^r f_Y(t)\,{\rm d}t,
	\label{eq:vr}
\end{equation}
and
\begin{equation}
	m_a(y)
	\coloneqq
	\int_{-\infty}^{\infty}
	\left[
	a\min\{y,t\}+(1-a)\max\{y,t\}
	\right]f_Y(t)\,{\rm d}t.
	\label{eq:ma}
\end{equation}

By construction, for  $Y_1\sim{\rm Gumbel}(0,1)$,  $Y_2\sim{\rm Gumbel}(0,1)$, independent, $Y_1-Y_2\sim\operatorname{Logistic}(0,1)$ (Lemma \ref{1-lemma}) and
\begin{align}\label{exp-imp}
	\begin{array}{lll}
\E[v_r(Y)]
&= \displaystyle
\E[|Y_1-Y_2|^r]
\stackrel{\eqref{eq:logisticmoment}}{=}
2\Gamma(r+1)\sum_{j=1}^{\infty}
\frac{(-1)^{j-1}}{j^r},
\\[0,5cm]
\E[m_a(Y)]
&= \displaystyle
a\E[\min\{Y_1,Y_2\}]
+
(1-a)\E[\max\{Y_1,Y_2\}]
\stackrel{\eqref{exp-Y}}{=}
(1-2a)\log(2)-\gamma.
	\end{array}
\end{align}

Define the first-order influence functions
\begin{equation}
	\phi_r(y)
	\coloneqq 
	-\frac{2}{2r\Gamma(r+1)\sum_{j=1}^{\infty}
		\frac{(-1)^{j-1}}{j^r}} \, 
	\left\{v_r(y)-2\Gamma(r+1)\sum_{j=1}^{\infty}
	\frac{(-1)^{j-1}}{j^r}\right\},
	\label{eq:phir}
\end{equation}
and
\begin{equation}
	\psi_{a,r}(y)
	\coloneqq 
	2\{m_a(y)-[(1-2a)\log(2)-\gamma]\}
	+
	[(1-2a)\log(2)-\gamma]\phi_r(y).
	\label{eq:psiar}
\end{equation}

The following theorem provides the joint
asymptotic distribution.

\begin{theorem}
	For every fixed $r>0$ and $a\in\R$,
	\[
	\sqrt n
	\begin{pmatrix}
		\log(\widehat{k}_r)-\log(k)
		\\[0,2cm]
		\log(\widehat{\lambda}_{a,r})-\log(\lambda)
	\end{pmatrix}
	\xrightarrow{d}
	N_2\left(
	\begin{pmatrix}0\\[0,2cm] 0\end{pmatrix},
	\boldsymbol{\Sigma}_{a,r}
	\right),
	\]
	where
	\begin{equation}\label{eq:Sigma}
		\boldsymbol{\Sigma}_{a,r}
		=
		\begin{pmatrix}
			\E[\phi_r(Y)^2] & \E[\phi_r(Y)\psi_{a,r}(Y)]/k
			\\[0,2cm]
			\E[\phi_r(Y)\psi_{a,r}(Y)]/k& \E[\psi_{a,r}(Y)^2]/k^2
		\end{pmatrix},
		\quad 
		Y\sim{\rm Gumbel}(0,1).
	\end{equation}
\end{theorem}
\begin{proof}
	{
	Consider the U-statistic kernel
	$
	h(y_1,y_2)=|y_1-y_2|^r.
	$
	Let $Y_1$ and $Y_2$ be independent random variables
	with standard Gumbel distribution \eqref{eq:gumbelcdf}.
	The Hoeffding decomposition expresses the kernel
	as the sum of its expectation, two first-order
	projections, and a second-order projection:
	\[
	\begin{aligned}
		h(y_1,y_2)
		=
		\mathbb{E}[|Y_1-Y_2|^r]+h_1(y_1)+h_1(y_2)
		+h_2(y_1,y_2),
	\end{aligned}
	\]
	where the first-order projection is
	\begin{align}\label{eq:vr11}
	h_1(y)
	=
	\mathbb{E}[|y-Y_2|^r]-\mathbb{E}[|Y_1-Y_2|^r]
	\stackrel{\eqref{eq:vr}}{=}
	v_r(y)-\mathbb{E}[|Y_1-Y_2|^r],
	\end{align}
	and the second-order projection is
	\[
	h_2(y_1,y_2)
	=
	|y_1-y_2|^r
	-
	\mathbb{E}[|Y_1-Y_2|^r]
	-
	h_1(y_1)
	-
	h_1(y_2).
	\]
	By construction,
	$
	\mathbb{E}[h_1(Y_1)]=0
	$
	and
	$
	\mathbb{E}[h_2(y,Y_2)]=0
	\ \text{for every }y.
	$
	Substituting this decomposition into the
	definition of the U-statistic gives
	\[
	\begin{aligned}
		k^r U_{n,r}
		=
		\mathbb{E}[|Y_1-Y_2|^r]
		+
		\frac{2}{n}\sum_{i=1}^n h_1(Y_i)
		+
		\frac{1}{\binom{n}{2}}
		\sum_{1\leqslant i<j\leqslant n}h_2(Y_i,Y_j).
	\end{aligned}
	\]
	The last term is the contribution of the
	second-order projection. Since this projection
	has conditional expectation zero, its variance
	is of order $O(n^{-2})$. Consequently, it is
	$o_p(n^{-1/2})$. Therefore, since $Y_1-Y_2\sim\operatorname{Logistic}(0,1)$ (Lemma \ref{1-lemma}),
	\begin{align} \label{dec-Uk}
		k^r U_{n,r}
		&\stackrel{\eqref{eq:vr11}}{=}
		\mathbb{E}[|Y_1-Y_2|^r]
		+
		\frac{2}{n}\sum_{i=1}^n
		\left\{v_r(Y_i)-\mathbb{E}[|Y_1-Y_2|^r]\right\}
		+
		R_n
		\nonumber
		\\[0,2cm]
		&\stackrel{\eqref{eq:logisticmoment}}{=}
		2\Gamma(r+1)\sum_{j=1}^{\infty}
		\frac{(-1)^{j-1}}{j^r}
		+
		\frac{2}{n}\sum_{i=1}^n
		\left\{
		v_r(Y_i)-2\Gamma(r+1)\sum_{j=1}^{\infty}
		\frac{(-1)^{j-1}}{j^r}
		\right\}
		+R_n
		\nonumber
				\\[0,2cm]
		&\stackrel{\eqref{eq:phir}}{=}
		2\Gamma(r+1)\sum_{j=1}^{\infty}
		\frac{(-1)^{j-1}}{j^r}
		\left[
		1
		-
		{\frac{r}{n}} % ORIGINAL: \frac{1}{n}
		\sum_{i=1}^n
		\phi_r(Y_i)
		\right]
		+
		R_n,
	\end{align}
	where $R_n=o_p(n^{-1/2})$, that is,
	$
	\sqrt{n}\,R_n\xrightarrow{p}0.
	$
	}
	
%	We first consider the U-statistic kernel
%	$
%	|y_1-y_2|^r,
%	$
%	$r>0$.
%	Since $Y_1-Y_2\sim\operatorname{Logistic}(0,1)$ (Lemma \ref{1-lemma}),
%	its first-order Hoeffding projection is
%	\[
%	\begin{aligned}
%	\mathbb{E}[|y-Y_2|^r]
%		-\mathbb{E}[|Y_1-Y_2|^r]
%		\stackrel{\eqref{eq:logisticmoment}}{=}
%		v_r(y)-2\Gamma(r+1)\sum_{j=1}^{\infty}
%		\frac{(-1)^{j-1}}{j^r}.
%	\end{aligned}
%	\]
%
%By the Hoeffding decomposition, the first-order
%expansion of the corresponding U-statistic is
%obtained by averaging twice this projection over
%the sample observations. More precisely,
%\[
%\begin{aligned}
%	k^r U_{n,r}
%	=
%	2\Gamma(r+1)\sum_{j=1}^{\infty}
%	\frac{(-1)^{j-1}}{j^r}
%	+
%	\frac{2}{n}\sum_{i=1}^n
%	\left\{
%	v_r(Y_i)-2\Gamma(r+1)\sum_{j=1}^{\infty}
%	\frac{(-1)^{j-1}}{j^r}
%	\right\}
%	+R_n,
%\end{aligned}
%\]
%where $R_n=o_p(n^{-1/2})$, meaning that
%$
%\sqrt{n}\,R_n\xrightarrow{p}0.
%$

To obtain the corresponding expansion for
$\log(\widehat{k}_r)$, we apply the delta method
to the function
\[
g(u)=\frac{1}{r}
\left\{
\log\left(2\Gamma(r+1)\sum_{j=1}^{\infty}
\frac{(-1)^{j-1}}{j^r}\right)
-\log(u)
\right\}.
\]
Since
\begin{align}\label{der-1}
g'\left(2\Gamma(r+1)\sum_{j=1}^{\infty}
\frac{(-1)^{j-1}}{j^r}\right)
=
-\frac{1}{2 r\Gamma(r+1)\sum_{j=1}^{\infty}
	\frac{(-1)^{j-1}}{j^r}}
\end{align}
and
$
g(2\Gamma(r+1)\sum_{j=1}^{\infty}
{(-1)^{j-1}}/{j^r})
=
0
$
a first-order Taylor expansion around
$2\Gamma(r+1)\sum_{j=1}^{\infty}
{(-1)^{j-1}}/{j^r}$, 
%together with the
%Hoeffding decomposition above, 
yields
\begin{align*}
	g(k^rU_{n,r})
	=
	g'\left(2\Gamma(r+1)\sum_{j=1}^{\infty}
	\frac{(-1)^{j-1}}{j^r}\right)
	\left[k^rU_{n,r}-2\Gamma(r+1)\sum_{j=1}^{\infty}
	\frac{(-1)^{j-1}}{j^r}\right]
	+
	o_p(n^{-1/2}).
\end{align*}
Using the identity $g(k^rU_{n,r})=\log(\widehat{k}_r)-\log(k)$ together with \eqref{der-1}, the above identity can be written as
\begin{align}\label{comb-1}
	\log(\widehat{k}_r)-\log(k)
	=
	-\frac{
		k^rU_{n,r}-
		2\Gamma(r+1)\sum_{j=1}^{\infty}
		\frac{(-1)^{j-1}}{j^r}
	}{
		2r\Gamma(r+1)\sum_{j=1}^{\infty}
		\frac{(-1)^{j-1}}{j^r}
	}
	+o_p(n^{-1/2})
	=
	\frac{1}{n}\sum_{i=1}^n\phi_r(Y_i)
	+o_p(n^{-1/2}),
\end{align}
where the last equality follows from the Hoeffding decomposition in \eqref{dec-Uk}.

	This representation identifies $\phi_r(Y_i)$
	as the first-order contribution of the $i$th
	observation to the estimation error of
	$\log(k)$. Such a function is called the
	influence function (IF). In particular, it
	describes how a small perturbation of the
	underlying distribution at a given observation
	affects the estimator, to first order.
	
	\smallskip 
	We next consider the U-statistic $V_{n,a}$ given in
	\eqref{eq:Vnr}, with kernel
	$
	h_a(x_1,x_2)
	=
	a\log(\min\{x_1,x_2\})
	+(1-a)\log(\max\{x_1,x_2\}).
	$
%	Using $\log(X_i)=\log(\lambda)+Y_i/k$ and the
%	fact that the logarithm is an increasing function,
%	we can write
%	\[
%	\begin{aligned}
%		h_a(X_1,X_2)
%		&=
%		\log(\lambda)
%		+
%		\frac{1}{k}
%		\left[
%		a\min\{Y_1,Y_2\}
%		+(1-a)\max\{Y_1,Y_2\}
%		\right].
%	\end{aligned}
%	\]
%	
%	Let $Y_1$ and $Y_2$ be independent random variables
%	with standard Gumbel distribution \eqref{eq:gumbelcdf}.
%	The expectation of the kernel is
%	\[
%	\begin{aligned}
%		\mathbb{E}[h_a(X_1,X_2)]
%		&=
%		\log(\lambda)
%		+\frac{1}{k} \, 
%		\mathbb{E}\left[
%		a\min\{Y_1,Y_2\}
%		+(1-a)\max\{Y_1,Y_2\}
%		\right]\\
%		&\stackrel{\eqref{eq:ma}}{=}
%		\log(\lambda)
%+
%{\mathbb{E}[m_a(Y)]\over k}
%		\stackrel{\eqref{exp-imp}}{=}
%		\log(\lambda)
%		+
%		\frac{(1-2a)\log(2)-\gamma}{k},
%	\end{aligned}
%	\]
%	where $Y\sim{\rm Gumbel}(0,1)$.
%	
	The first-order Hoeffding projection is therefore
	\begin{align}\label{proj-1}
		h_{a,1}(y)
		=
		\mathbb{E}[h_a(X_1,X_2)\mid Y_1=y]
		-\mathbb{E}[h_a(X_1,X_2)]
		=
		\frac{1}{k}
		\left\{
		m_a(y)-[(1-2a)\log(2)-\gamma]
		\right\}.
	\end{align}

	The second-order projection is given by
	\[
	\begin{aligned}
		h_{a,2}(y_1,y_2)
		=
		h_a(x_1,x_2)
		-
		\mathbb{E}[h_a(X_1,X_2)]
		-
		h_{a,1}(y_1)
		-
		h_{a,1}(y_2),
	\end{aligned}
	\]
	where $x_i=\lambda\exp(y_i/k)$, for $i=1,2$.
	By construction,
	$
	\mathbb{E}[h_{a,1}(Y_1)]=0
	\ \text{and} \
	\mathbb{E}[h_{a,2}(y,Y_2)]=0
	\ 
	\text{for every }y.
	$
	
	Substituting this decomposition into the definition
	of the U-statistic gives
	\[
	\begin{aligned}
		V_{n,a}
		=
		\mathbb{E}[h_a(X_1,X_2)]
		+\frac{2}{n}\sum_{i=1}^n h_{a,1}(Y_i)
		+
		\frac{1}{\binom{n}{2}}
		\sum_{1\leqslant i<j\leqslant n}
		h_{a,2}(Y_i,Y_j).
	\end{aligned}
	\]
	The last term is the contribution of the
	second-order projection. Since this projection
	has conditional expectation zero and a finite
	second moment, its variance is of order
	$O(n^{-2})$. Consequently, it is
	$o_p(n^{-1/2})$. Hence, by \eqref{eq:Cexpectation} and \eqref{proj-1},
	\begin{align}\label{dec-Vna}
		V_{n,a}
		=
		\log(\lambda)
		+
		\frac{(1-2a)\log(2)-\gamma}{k}
		+
		\frac{2}{nk}\sum_{i=1}^n
		\left\{
		m_a(Y_i)-[(1-2a)\log(2)-\gamma]
		\right\}
		+R_{n,a},
	\end{align}
	where $R_{n,a}=o_p(n^{-1/2})$, that is,
	$\sqrt{n}\,R_{n,a}\xrightarrow{p}0$.
	
	To obtain the corresponding expansion for
	$\log(\widehat{\lambda}_{a,r})$, we use the
	Hoeffding decomposition in \eqref{dec-Vna}.
	Since
	\[
	\log(\widehat{\lambda}_{a,r})
	=
	V_{n,a}
	-\frac{(1-2a)\log(2)-\gamma}{\widehat{k}_r},
	\]
	subtracting $\log(\lambda)$ and using
	\eqref{dec-Vna} yield
	\begin{align}\label{exp-V}
		\log(\widehat{\lambda}_{a,r})-\log(\lambda)
		&=
		\frac{2}{nk}\sum_{i=1}^n
		\left\{
		m_a(Y_i)-[(1-2a)\log(2)-\gamma]
		\right\}
		\nonumber
		\\[0,2cm]
		&+
		\frac{(1-2a)\log(2)-\gamma}{k}
		-\frac{(1-2a)\log(2)-\gamma}{\widehat{k}_r}
		+R_{n,a}
		\nonumber
		\\[0,2cm]
		&=
		\frac{2}{nk}\sum_{i=1}^n
		\left\{
		m_a(Y_i)-[(1-2a)\log(2)-\gamma]
		\right\}
		\nonumber
		\\[0,2cm]
		&
        +
		[(1-2a)\log(2)-\gamma]
		\left(
		\frac{1}{k}-\frac{1}{\widehat{k}_r}
		\right)
		+
		R_{n,a}.
	\end{align}
	
	Now, using
	\eqref{comb-1} and a first-order Taylor
	expansion of the exponential function around
	zero, we obtain
	\[
	\begin{aligned}
		\frac{1}{\widehat{k}_r}
		=
		\frac{1}{k}
		\exp\left\{
		-[\log(\widehat{k}_r)-\log(k)]
		\right\}
		=
		\frac{1}{k}
		-\frac{1}{nk}\sum_{i=1}^n\phi_r(Y_i)
		+o_p(n^{-1/2}).
	\end{aligned}
	\]
	Therefore,
	\[
	\frac{1}{k}-\frac{1}{\widehat{k}_r}
	=
	\frac{1}{nk}\sum_{i=1}^n\phi_r(Y_i)
	+o_p(n^{-1/2}).
	\]
	
	Substituting this expression into the expansion \eqref{exp-V}
	gives
	\begin{align}\label{comb-2}
		\log(\widehat{\lambda}_{a,r})-\log(\lambda)
		&=
		\frac{1}{nk}\sum_{i=1}^n
		\Big[
		2\{m_a(Y_i)-[(1-2a)\log(2)-\gamma]\}
		\nonumber
		\\[0,2cm]
		&
		+
		[(1-2a)\log(2)-\gamma]\phi_r(Y_i)
		\Big]
		+o_p(n^{-1/2})
		\nonumber
		\\[0,2cm]
		&=
		\frac{1}{nk}\sum_{i=1}^n
		\psi_{a,r}(Y_i)+o_p(n^{-1/2}),
	\end{align}
	where $\psi_{a,r}$ is defined in \eqref{eq:psiar}.
	Thus, the influence function for
	$\log(\widehat{\lambda}_{a,r})$ is
	$\psi_{a,r}(y)/k$.

By construction,
$
\mathbb{E}[\phi_r(Y)]=0,
\
\mathbb{E}[\psi_{a,r}(Y)]=0,
$
and
$
\mathbb{E}[\phi_r(Y)^2]<\infty,
\
\mathbb{E}[\psi_{a,r}(Y)^2]<\infty,
$
for every $r>0$ and $a\in\mathbb{R}$, since the
Gumbel distribution has finite moments of every
positive order. Combining the two asymptotic
linear representations \eqref{comb-1} and \eqref{comb-2}, we obtain
\[
\begin{aligned}
	\sqrt{n}
	\begin{pmatrix}
		\log(\widehat{k}_r)-\log(k)
		\\[0,2cm]
		\log(\widehat{\lambda}_{a,r})-\log(\lambda)
	\end{pmatrix}
	=
	\frac{1}{\sqrt{n}}\sum_{i=1}^n
	\begin{pmatrix}
		\phi_r(Y_i)
		\\[0,2cm]
		\psi_{a,r}(Y_i)/k
	\end{pmatrix}
	+o_p(1).
\end{aligned}
\]
Since $Y_1,\ldots,Y_n$ are i.i.d., the vectors
in the summation are also i.i.d. Therefore, the
multivariate central limit theorem implies that
the normalized sum converges in distribution to
a bivariate normal distribution with zero mean
and covariance matrix $\boldsymbol{\Sigma}_{a,r}$ given in \eqref{eq:Sigma}.
The asserted result then follows from
Slutsky's theorem.
\end{proof}

\begin{remark}
	The asymptotic covariance matrix $\boldsymbol{\Sigma}_{a,r}$ in
	\eqref{eq:Sigma} is expressed for the logarithms
	of the estimators. This parametrization is
	convenient because the shape and scale parameters
	are positive and because the covariance structure
	has a simple dependence on the unknown parameters.
\end{remark}

\section{Joint asymptotic optimization}
\label{Joint asymptotic optimization}

The tuning parameters $a$ and $r$ determine the
asymptotic covariance matrix $\boldsymbol{\Sigma}_{a,r}$ of the proposed
estimators. We consider the determinant of this
matrix as a joint measure of asymptotic dispersion.

Let
\begin{align}\label{def-D}
	D(a,r)
\coloneqq
	\E[\phi_r(Y)^2] \E[\psi_{a,r}(Y)^2] -\E[\phi_r(Y)\psi_{a,r}(Y)]^2.
\end{align}
From \eqref{eq:Sigma},
\begin{equation*}
	\det\{\boldsymbol{\Sigma}_{a,r}\}
	=
	\frac{1}{k^2} \, D(a,r).
%	\label{eq:detSigma}
\end{equation*}

Since $k$ does not depend on the tuning parameters,
the optimization problem reduces to
\begin{equation}
	\min_{a\in\R,\;r>0}
D(a,r).
	\label{eq:optimization}
\end{equation}

\subsection{Analytical reduction with respect to $a$}

%Let
%$$
%m_0(y)
%\coloneqq 
%\E[\max\{y,Y\}]
%\quad 
%\text{and}
%\quad 
%d(y)
%\coloneqq 
%\E[\min\{y,Y\}]-\E[\max\{y,Y\}].
%$$
Since $m_a(y)$ in \eqref{eq:ma} can be expressed as
$
m_a(y)=\E[\max\{y,Y\}]+a\{\E[\min\{y,Y\}]-\E[\max\{y,Y\}]\},
$ with $Y\sim{\rm Gumbel}(0,1)$, $\psi_{a,r}(y)$ in \eqref{eq:psiar}
%and
%$
%(1-2a)\log(2)-\gamma=\theta_0+a\theta_1,
%$
%where
%$
%\theta_0=-\gamma+\log(2),
%\
%\theta_1=-2\log(2),
%$
can be written as
\[
\psi_{a,r}(y)=P_r(y)+aQ_r(y),
\]
with
\begin{align*}
	P_r(y)
	&\coloneqq
	2\{\E[\max\{y,Y\}]-\log(2)+\gamma\}
	+[\log(2)-\gamma]\phi_r(y),
	\\[0,2cm]
	Q_r(y)
	&\coloneqq
	2\left\{ \E[\min\{y,Y\}]-\E[\max\{y,Y\}]+2\log(2)\right\}
	-2\log(2)\phi_r(y).
\end{align*}

%Define
%\begin{align*}
%	&V_0(r)\coloneqq\E[P_r(Y)^2],\quad  
%	V_1(r)\coloneqq\E[P_r(Y)Q_r(Y)],\quad
%	V_2(r)\coloneqq\E[Q_r(Y)^2],
%	\\[0,2cm]
%	&J_0(r)\coloneqq\E[\phi_r(Y)P_r(Y)],\quad
%	J_1(r)\coloneqq\E[\phi_r(Y)Q_r(Y)].
%\end{align*}

It follows that
\begin{align*}
	&\E[\psi_{a,r}(Y)^2] =\E[P_r(Y)^2]+2a\E[P_r(Y)Q_r(Y)]+a^2\E[Q_r(Y)^2],
	\\[0,2cm]
	&\E[\phi_r(Y)\psi_{a,r}(Y)]=\E[\phi_r(Y)P_r(Y)]+a \E[\phi_r(Y)Q_r(Y)].
\end{align*}

Consequently, $D(a,r)$ in \eqref{def-D} can be written as
\begin{equation}
	D(a,r)
	=
	A(r)a^2+2B(r)a+C(r),
	\label{eq:Dquadratic}
\end{equation}
where
\begin{align*}
	A(r)&\coloneqq\E[\phi_r(Y)^2] \E[Q_r(Y)^2]-\E[\phi_r(Y)Q_r(Y)]^2,
	\\[0,2cm]
	B(r)&\coloneqq\E[\phi_r(Y)^2] \E[P_r(Y)Q_r(Y)]-\E[\phi_r(Y)P_r(Y)] \E[\phi_r(Y)Q_r(Y)],
	\\[0,2cm]
	C(r)&\coloneqq\E[\phi_r(Y)^2] \E[P_r(Y)^2]-\E[\phi_r(Y)P_r(Y)]^2.
\end{align*}

\begin{proposition}
	Suppose that $A(r)>0$ for a given $r>0$.
	Then the unique minimizer of $D(a,r)$ in \eqref{def-D} over
	$a\in\mathbb{R}$ is
	\begin{equation}
			a^*(r)=-\frac{B(r)}{A(r)}.
		\label{eq:astar}
	\end{equation}
	provided $A(r)>0$.
	The corresponding minimized criterion is
	\begin{equation}
			D^*(r)=C(r)-\frac{B(r)^2}{A(r)}.
		\label{eq:Dstar}
	\end{equation}
\end{proposition}

\begin{proof}
	By \eqref{eq:Dquadratic}, $D(a,r)$ is a
	quadratic function of $a$. Differentiating
	with respect to $a$ gives
	\[
	\frac{\partial D(a,r)}{\partial a} 
	=
	2A(r)a+2B(r).
	\]
	The unique stationary point is therefore
	$
	a^*(r)=-B(r)/A(r).
	$
	Since $A(r)>0$, this point is the unique
	global minimizer for fixed $r$.
	Substitution into \eqref{eq:Dquadratic}
	gives \eqref{eq:Dstar}.
\end{proof}

{
\begin{remark}\label{rem:A1}
The hypothesis $A(r)>0$ fails at $r=1$. Actually,
$$\E[\min\{y,Y\}]-\E[\max\{y,Y\}]=-v_1(y)$$ and, by \eqref{eq:phir},
$$\phi_1(y)=-\{v_1(y)-2\log(2)\}/\log(2),$$
hence $Q_1\equiv 0$.
Therefore, $A(1)=B(1)=0$ and $D(a,1)=C(1)$ for every $a\in\mathbb{R}$,
in agreement with $\widehat{\lambda}_{a,1}$ does not
depending on $a$ (Section~\ref{Model and proposed estimators}).
Moreover, $a^*(r)=-B(r)/A(r)$ is an indeterminate form $0/0$ as
$r\to1$, hence it must be evaluated with care in a neighbourhood
of $r=1$.
\end{remark}}

%\begin{remark}
%	The Cauchy--Schwarz inequality implies
%	$
%	A(r)\geqslant 0.
%	$
%	Strict positivity must be verified for the
%	values of $r$ under consideration. It should
%	not be assumed without checking the possible
%	linear dependence of the corresponding
%	influence functions.
%\end{remark}

% %------------------------------------------------
% \subsection{One-dimensional optimization}
% %------------------------------------------------

\begin{remark}
    The original problem \eqref{eq:optimization}
can now be reduced to
\begin{equation}
		r^* \in
		\argmin_{r>0}D^*(r),
	\label{eq:rstar}
\end{equation}
and $a^*=a^*(r^*).$ If the minimum is attained at an interior
point and $D^*$ is differentiable, then
$
{\rm d}D^*(r^*)/{\rm d} r =0.
$

The existence, uniqueness, and numerical
location of the minimizer require further
investigation. In particular, the behavior
of $D^*(r)$ as $r\to0^+$ and $r\to\infty$
must be examined before a global optimum
can be claimed.
\end{remark}

\section{Numerical investigation}
\label{Numerical investigation}

In this section, we present a Monte Carlo simulation study to evaluate the derived estimators $\widehat{k}_r$ and $\widehat{\lambda}_{a,r}$, described in \eqref{eq:kestimator} and \eqref{eq:lambdaestimator}.
The simulation study was programmed in the software R.

We consider the following fixed parameters for the simulations:
\begin{itemize}
\item true population parameters $k\in\{0.5, 1, 2.5, 5\}$ and $\lambda=1$.
    \item $M=1,000$ Monte Carlo replications of samples $\{X_1, \cdots, X_n\}$ sampled from $\operatorname{Weibull}(k, \lambda)$;
    \item  sample size $n\in\{10, 20, 50, 100, 300\}$;
\end{itemize}

To study the Monte Carlo to evaluate the estimators, we compute the relative bias (RB) and root mean squared error (RMSE) of each estimated value compared to its true population parameter. 
We also compute the coverage probability (CP) of the estimates, which measures the proportion of times the confidence intervals contain the true value of the parameter.
%The results were compared with the modified MLE (MMLE) and corrected MLE (CMLE), which were proposed by \citep{kundu2009estimation} and \citep{cheng1987corrected}, respectively.

\subsection{Selecting $(r,a)$ from a fixed grid}

In this subsection, we evaluate the impact of selecting $(r,a)$ from a fixed grid. 

\begin{figure}[H]
    \centering
\includegraphics[width=1.0\linewidth]{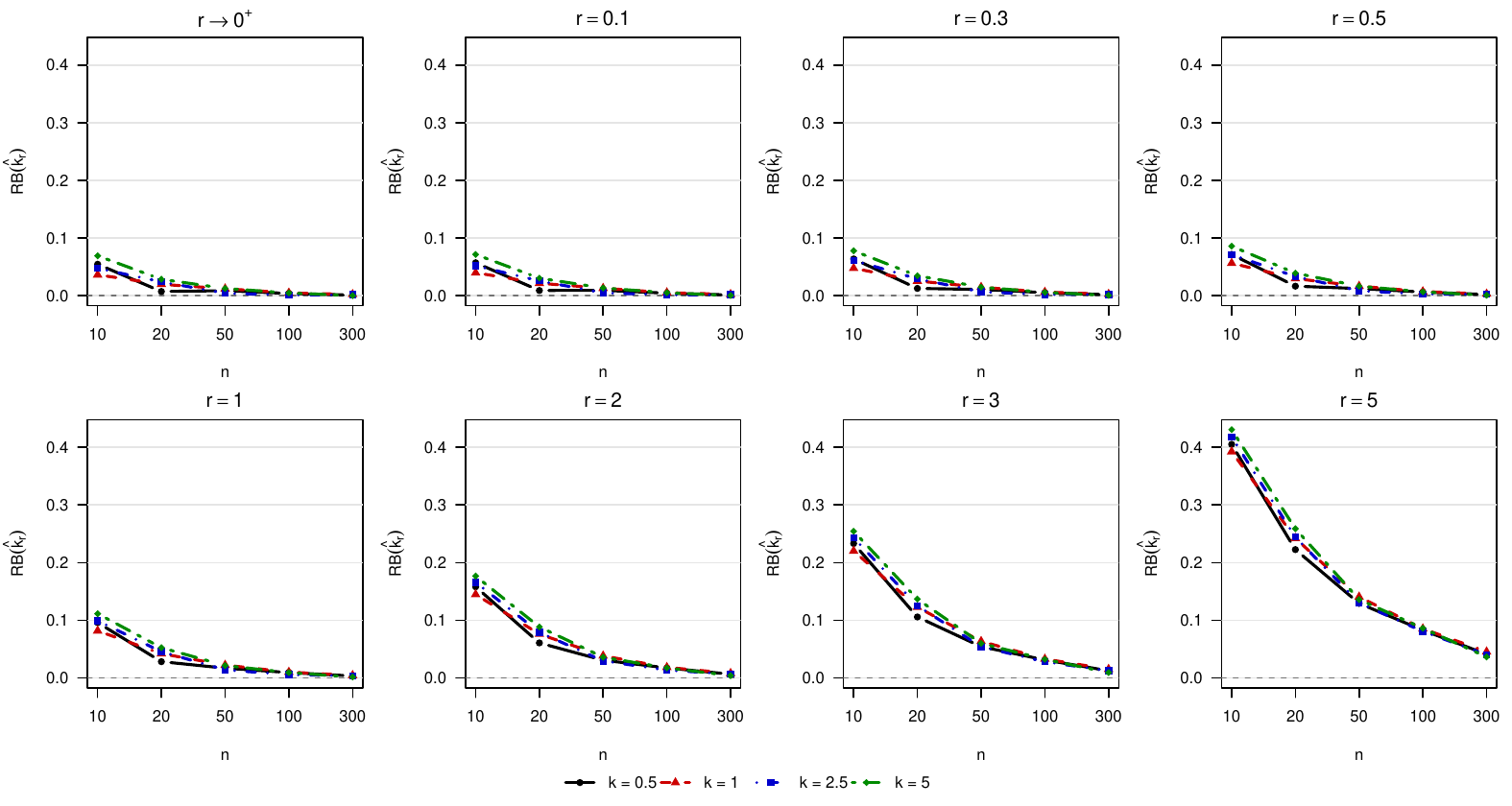}
    \caption{ Monte Carlo results for the RB decay of $\widehat{k}_r$ estimates of $k\in\{0.5, 1, 2.5, 5\}$, with increasing sample size $n = 10, 20, 50, 100, 300$ and varying $r \in \{0^+, 0.1, 0.3, 0.5, 1, 2, 3, 5\}$. }
    \label{fig:sim11}
\end{figure}

\begin{figure}[H]
    \centering
\includegraphics[width=1.0\linewidth]{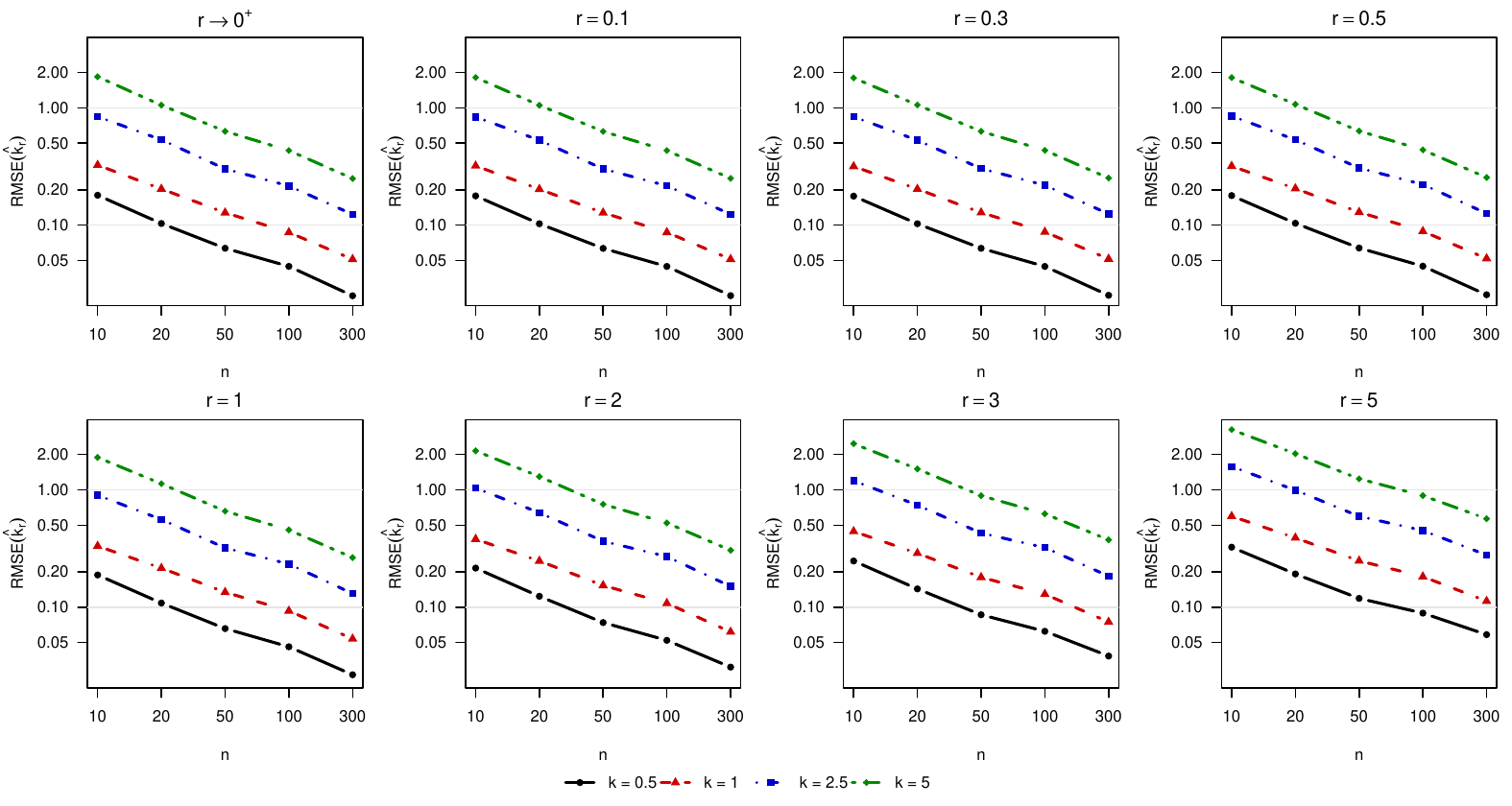}
    \caption{Monte Carlo results for the RMSE decay of $\widehat{k}_r$ estimates of $k\in\{0.5, 1, 2.5, 5\}$, with increasing sample size $n = 10, 20, 50, 100, 300$ and varying $r \in \{0^+, 0.1, 0.3, 0.5, 1, 2, 3, 5\}$.}
    \label{fig:sim12}
\end{figure}

\begin{figure}[H]
    \centering
\includegraphics[width=1.0\linewidth]{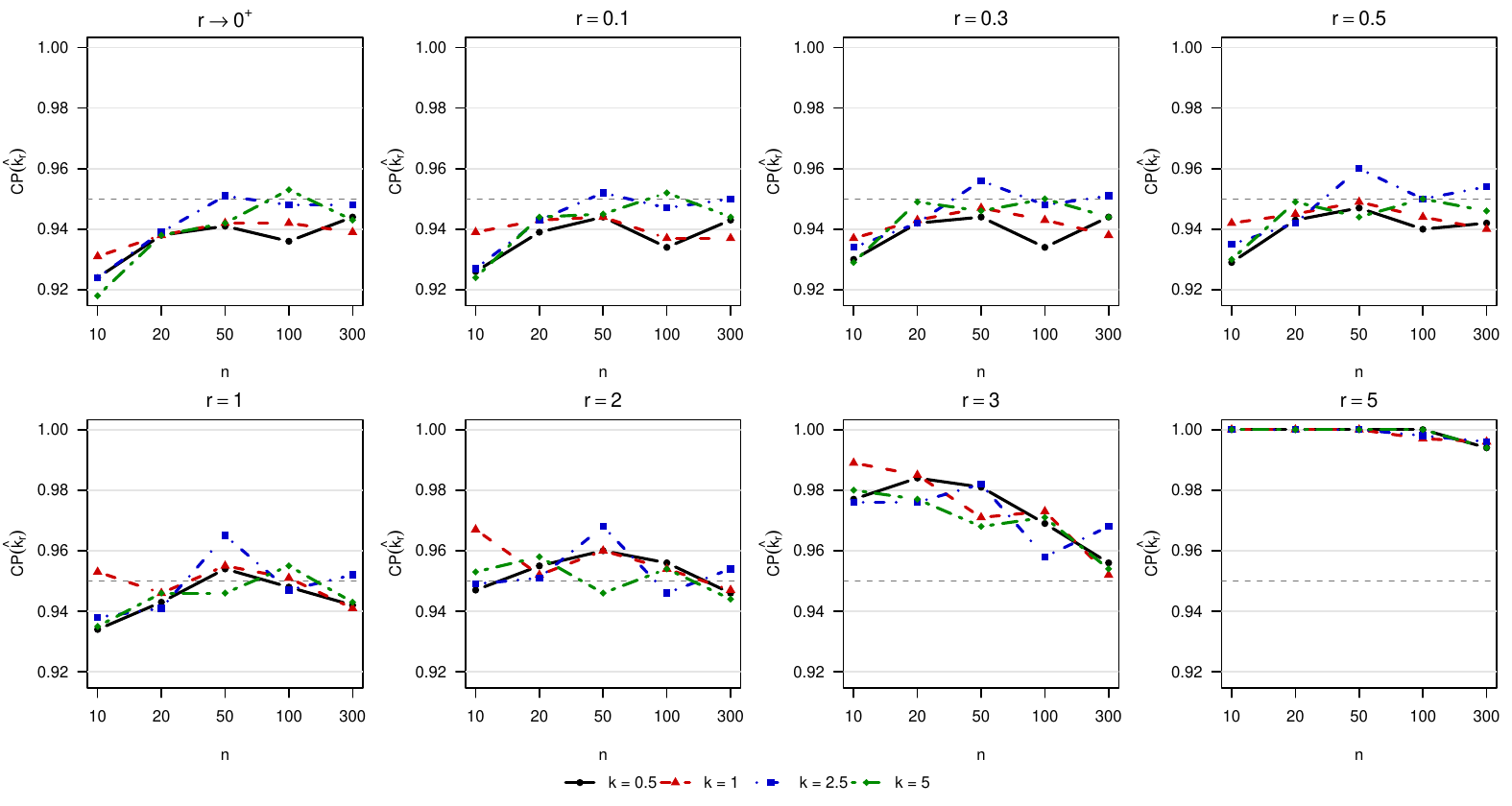}
    \caption{Monte Carlo results for the CP of $\widehat{k}_r$ estimates of $k\in\{0.5, 1, 2.5, 5\}$, with increasing sample size $n = 10, 20, 50, 100, 300$ and varying $r \in \{0^+, 0.1, 0.3, 0.5, 1, 2, 3, 5\}$.}
    \label{fig:sim13}
\end{figure}

First, considering the grid $r \in \{0^+, 0.1, 0.3, 0.5, 1, 2, 3, 5\}$, Figures \ref{fig:sim11}--\ref{fig:sim13} display the RB, RMSE, and CP of the estimator $\widehat{k}_r$, respectively. The RB exhibits its best performance as $r \to 0^+$, while RB values exceeding $0.1$ are observed for small sample sizes when $r > 2$. Overall, both RB and RMSE converge to zero as $n$ increases. For $r = 5$, larger sample sizes are required for the RMSE to reach lower levels. Regarding CP, the coverage probabilities approach the nominal level of $0.95$ when $r < 3$ and $n \approx 300$. However, for $r \geq 3$, the CP values tend to exceed $0.95$, indicating wider confidence intervals. 

\begin{figure}[H]
    \centering
\includegraphics[width=1.0\linewidth]{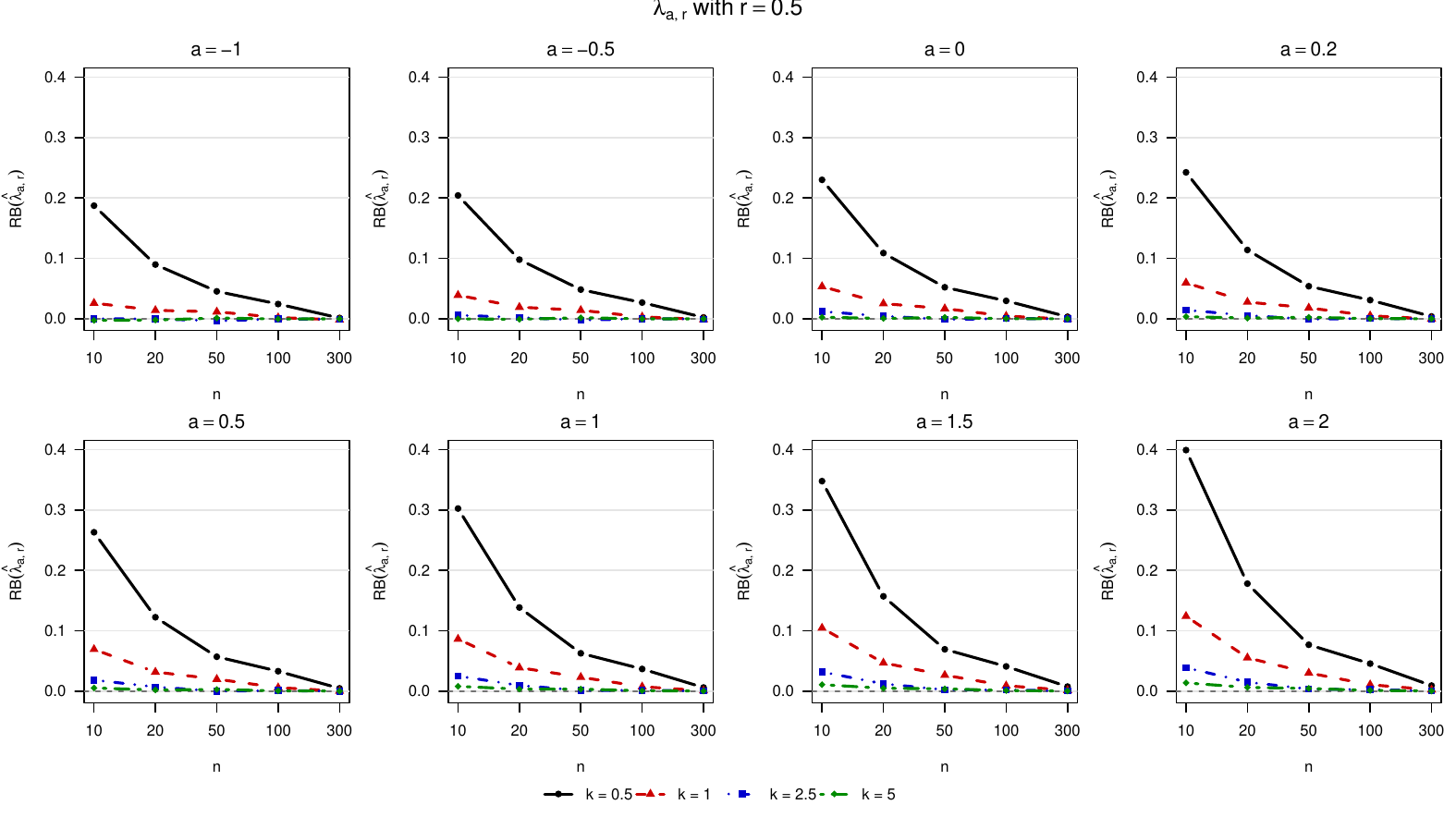}
    \caption{Monte Carlo results for the RB decay of $\widehat{\lambda}_{a,0.5}$ estimates of $\lambda=1$, with increasing sample size $n = 10, 20, 50, 100, 300$ and varying $a \in \{-1, -0.5, 0, 0.2, 0.5, 1.5, 2\}$.
}
    \label{fig:sim14}
\end{figure}

\begin{figure}[H]
    \centering
\includegraphics[width=1.0\linewidth]{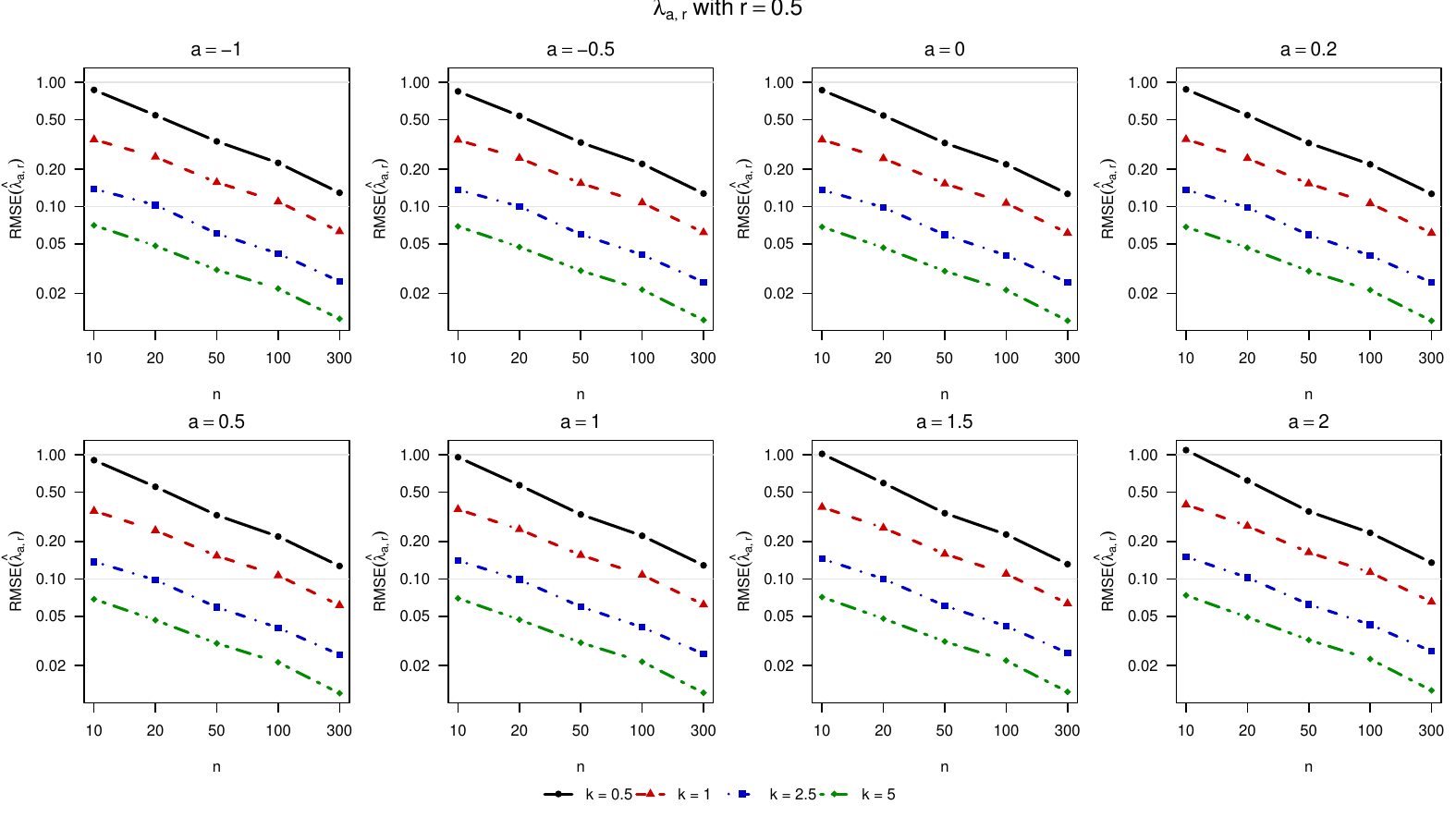}
    \caption{Monte Carlo results for the RMSE decay of $\widehat{\lambda}_{a,0.5}$ estimates of $\lambda=1$, with increasing sample size $n = 10, 20, 50, 100, 300$ and varying $a \in \{-1, -0.5, 0, 0.2, 0.5, 1.5, 2\}$.}
    \label{fig:sim15}
\end{figure}

\begin{figure}[H]
    \centering
\includegraphics[width=1.0\linewidth]{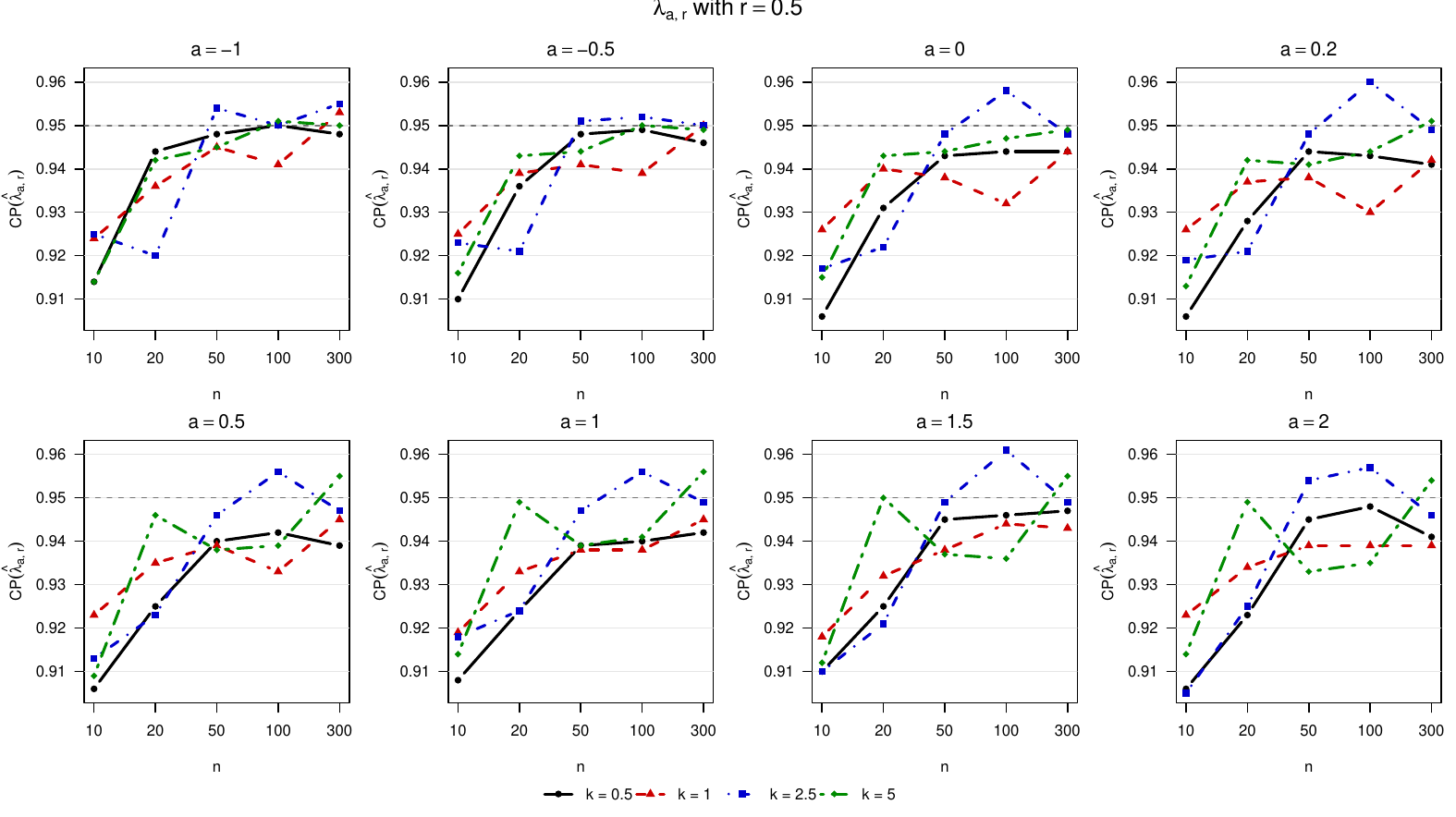}
    \caption{Monte Carlo results for the CP of $\widehat{\lambda}_{a,0.5}$ estimates of $\lambda=1$, with increasing sample size $n = 10, 20, 50, 100, 300$ and varying $a \in \{-1, -0.5, 0, 0.2, 0.5, 1.5, 2\}$.}
    \label{fig:sim16}
\end{figure}

Secondly, for a fixed value of $r = 0.5$, Figures \ref{fig:sim14}--\ref{fig:sim16} present the performance of the estimator $\widehat{\lambda}_{a,0.5}$ for $a \in \{-1, -0.5, 0, 0.2, 0.5, 1.5, 2\}$. In contrast to the results obtained for $\widehat{k}_r$, both RMSE and RB tend to be larger when $k$ is small. Nevertheless, as expected, RB and RMSE decrease as the sample size $n$ increases. For all values of $a$ considered, the CP remains close to the nominal level, with $0.94 \leq \mathrm{CP} \leq 0.96$.

\subsection{Comparison of two-step and joint parameter tuning under different optimizers}

In this subsection, in contrast to the previous one, we investigate the tuning of the parameters through optimization-based methods and compare the results with those obtained using the fixed reference table (see Appendix \ref{ap1}). We consider two optimization strategies: a two-step approach and a joint optimization approach.

Consider $M$ independent random samples from the standard $\operatorname{Weibull}(1,1)$ distribution. Let \[ \mathbf{W}^{(m)}=\left(W_1^{(m)}, \ldots, W_n^{(m)}\right)^\top, \qquad m=1,\ldots,M, \] denote the $m$-th random sample, and let $\widehat{k}_r^{(m)}$ be the corresponding estimate of $\kappa$. In the first step, we determine \[ \widetilde{r}_n = \arg\min_{r>0} n\,\mathrm{MSE}\!\left(\log\widehat{k}_r\right) \approx \arg\min_{r>0} \frac{n}{M}\sum_{m=1}^{M} \left(\log\widehat{k}^{(m)}_r\right)^2. \] Subsequently, we optimize
\begin{equation}\label{eq:a_tunning}
    \widetilde{a}_{n,r} \in\argmin_{a\in\mathbb{R}} n\,\mathrm{MSE}(\log\widehat{\lambda}_{a,r}) \approx \argmin_{a\in\mathbb{R}} \frac{n}{M}\sum_{m=1}^{M}\big(\log\widehat{\lambda}^{(m)}_{a,r}\big)^2.
\end{equation}

Note that the right-hand side of \eqref{eq:a_tunning} can be rewritten as
$$\frac{n}{M}\sum_{m=1}^{M}\big(\log\widehat{\lambda}^{(m)}_{a,r}\big)^2
=\frac{n}{M}\sum_{m=1}^{M}\left(\widehat{P}^{(m)}_r+a\,\left(\frac{2\log(2)}{\widehat{k}^{(m)}_r}-U^{(m)}_{n,1}\right)\right)^2,$$
where $\widehat{P}^{(m)}_r
	=\frac{1}{n}\sum_{i=1}^{n}\log W^{(m)}_i+\tfrac{1}{2}U^{(m)}_{n,1}-\frac{\log(2)-\gamma}{\widehat{k}^{(m)}_r}$, which gives us the closed-form expression 
$$	\widetilde{a}_{n,r}=-\frac{\sum_{m=1}^{M}\widehat{P}^{(m)}_r\left(\frac{2\log(2)}{\widehat{k}^{(m)}_r}-U^{(m)}_{n,1}\right)}
	{\sum_{m=1}^{M}\left(\frac{2\log(2)}{\widehat{k}^{(m)}_r}-U^{(m)}_{n,1}\right)^2}.$$

Based on $M=20\,000$ samples from the standard Weibull distribution, we constructed the reference table in Table~\ref{tab:lookup}, containing precomputed values of $(r,a)$. As established in Remark~\ref{remark:exp}, these values do not depend on $(k,\lambda)$ or on the observed data and can therefore be used universally.

\begin{remark}\label{remark:exp}
    Since $X=\lambda W^{1/k}$, with $W$ following the standard exponential
distribution, the errors $\log\widehat{k}_r-\log k$ and
$k(\log\widehat{\lambda}_{a,r}-\log\lambda)$ have the same distribution for every
$(k,\lambda)$. %, so that $n\,\mathrm{MSE}(\log\widehat{k}_r)$ depends only on $(n,r)$, and $n\,\mathrm{MSE}(\log\widehat{\lambda}_{a,r})$ depends on $(k,\lambda)$ only through the factor $k^{-2}$, which does not affect its minimizer.
Consequently, the selected pair $(r,a)$ depends only on the sample size, and not on
the parameters or on the observed data. %: it is computed once from simulated standard Weibull samples ($k=\lambda=1$), as in Table~\ref{tab:tuning}, and can be applied directly to real data sets of the same size, whose parameter values are unknown.
\end{remark}

\begin{remark}
 For the joint parameter tuning, we minimize $n^2$ times the determinant of the MSE matrix of $(\log\widehat{k}_r,\log\widehat{\lambda}_{a,r})$.
 
%  joint: $(r,a)$ minimizes $n^2\det\{\mathbf{M}_n(a,r)\}$, where
% $\mathbf{M}_n(a,r)=\E(\mathbf{e}\mathbf{e}^\top)$ is the $2\times 2$ mean squared error
% matrix of $\mathbf{e}=(\log\widehat{k}_r-\log k,\ \log\widehat{\lambda}_{a,r}-\log\lambda)^\top$,
% the finite-sample counterpart of $D(a,r)$ in \eqref{def-D};
\end{remark}

In Table~\ref{tab:tuning}, we compare the two-step and joint parameter-tuning approaches using the optimization methods grid search, Brent, BFGS, L-BFGS-B, and \texttt{nlminb}. For each method, the pair $(r_{\mathrm{opt}}, a_{\mathrm{opt}})$ obtained through optimization is compared with the pair $(\widetilde{r}_n,\widetilde{a}_{n,\widetilde{r}_n})$ provided by the reference table. Table~\ref{tab:tuning} reports the squared Euclidean distance \[ d^2=\left\|(r_{\mathrm{opt}},a_{\mathrm{opt}})-(\widetilde{r}_n,\widetilde{a}_{n,\widetilde{r}_n})\right\|^2. \] For two different parameter settings, the results indicate that the optimized values are close to those obtained from the reference table under both the two-step and joint optimization approaches. The computational time required by each optimization method is also reported.

Table~\ref{tab:tuning} also reports the percentage of runs that converged (Conv.), the percentage of runs in which the optimizer, using a single starting value, reached the minimum of the objective function within $0.01\%$ of the best value found (Single), and the percentage of sample sizes for which at least one of the 10 starting values attained this minimum (Multi).

\begin{table}[H]
\centering
\caption{Monte Carlo comparison of tuning-parameter optimization methods: distance to the reference table, computation time, and convergence rates.}
\label{tab:tuning}
%\scriptsize
\small
\setlength{\tabcolsep}{4pt}
\begin{tabular}{llcccccccc}
\toprule
& & \multicolumn{3}{c}{Runs (\%)} & \multicolumn{5}{c}{$n$} \\
\cmidrule(lr){3-5}\cmidrule(lr){6-10}
Parameters& Method & Conv. & Single & Multi & 10 & 20 & 50 & 100 & 300 \\
\midrule
% \multicolumn{10}{l}{\textit{(A) $k=\lambda=1$: $d^2$ (mean time per run in seconds)}} \\
$k=\lambda=1$ & \multicolumn{9}{l}{Two-step} \\
& \quad Grid search & -- & 60 & -- & 0.00134 & 0.01309 & 0.02348 & 0.00188 & 0.00002 \\
&  &  &  &  & (0.234) & (0.301) & (0.628) & (2.37) & (17.3) \\
& \quad Brent & -- & 100 & -- & 0.00115 & 0.01063 & 0.02225 & 0.00182 & 0.00002 \\
&  &  &  &  & (0.0680) & (0.0685) & (0.158) & (0.887) & (7.46) \\
& \quad BFGS & 94 & 86 & 100 & 0.00115 & 0.01063 & 0.02225 & 0.00118 & 0.00002 \\
&  &  &  &  & (0.406) & (0.381) & (0.533) & (2.16) & (16.5) \\
& \quad L-BFGS-B & 100 & 100 & 100 & 0.00115 & 0.01063 & 0.02225 & 0.00184 & 0.00002 \\
&  &  &  &  & (0.165) & (0.161) & (0.361) & (1.01) & (4.19) \\
& \quad \texttt{nlminb} & 100 & 100 & 100 & 0.00115 & 0.01063 & 0.02225 & 0.00182 & 0.00002 \\
&  &  &  &  & (0.122) & (0.137) & (0.259) & (0.779) & (3.23) \\
& \multicolumn{9}{l}{Joint} \\
& \quad Grid search & -- & 60 & -- & 0.00160 & 0.01570 & 0.02250 & 0.00170 & 0.00010 \\
&  &  &  &  & (0.239) & (0.439) & (0.537) & (2.61) & (17.8) \\
& \quad Profile (Brent) & -- & 100 & -- & 0.00004 & 0.01919 & 0.02642 & 0.00161 & 0.00001 \\
&  &  &  &  & (0.0445) & (0.101) & (0.138) & (1.03) & (8.47) \\
& \quad Nelder--Mead & 100 & 100 & 100 & 0.00004 & 0.01920 & 0.02642 & 0.00162 & 0.00001 \\
&  &  &  &  & (0.336) & (0.508) & (1.14) & (4.90) & (36.8) \\
& \quad BFGS & 100 & 90 & 100 & 0.00004 & 0.01919 & 0.02643 & 0.00076 & 0.00002 \\
&  &  &  &  & (0.420) & (0.566) & (1.07) & (6.36) & (61.6) \\
& \quad L-BFGS-B & 98 & 98 & 100 & 0.00004 & 0.01919 & 0.02643 & 0.00162 & 0.00001 \\
&  &  &  &  & (0.286) & (0.350) & (0.834) & (1.62) & (11.3) \\
& \quad \texttt{nlminb} & 100 & 100 & 100 & 0.00004 & 0.01919 & 0.02643 & 0.00162 & 0.00001 \\
&  &  &  &  & (0.229) & (0.267) & (0.710) & (1.39) & (9.82) \\
\midrule
%\multicolumn{9}{l}{\textit{(B) $k=1.5$, $\lambda=2$: $d^2$ (mean time per run in seconds)}} \\
$k=1.5$, $\lambda=2$ & \multicolumn{9}{l}{Two-step} \\
&\quad Grid search & -- & 60 & -- & 0.00134 & 0.01309 & 0.02348 & 0.00188 & 0.00002 \\
 &&  &  &  & (0.253) & (0.283) & (0.817) & (2.41) & (20.8) \\
&\quad Brent & -- & 100 & -- & 0.00115 & 0.01063 & 0.02225 & 0.00182 & 0.00002 \\
 &&  &  &  & (0.0615) & (0.0910) & (0.166) & (0.888) & (8.66) \\
&\quad BFGS & 94 & 86 & 100 & 0.00115 & 0.01063 & 0.02225 & 0.00118 & 0.00002 \\
 &&  &  &  & (0.425) & (0.401) & (0.569) & (2.06) & (19.1) \\
&\quad L-BFGS-B & 100 & 100 & 100 & 0.00115 & 0.01063 & 0.02225 & 0.00184 & 0.00002 \\
& &  &  &  & (0.156) & (0.180) & (0.420) & (0.898) & (4.75) \\
&\quad \texttt{nlminb} & 100 & 100 & 100 & 0.00115 & 0.01063 & 0.02225 & 0.00182 & 0.00002 \\
& &  &  &  & (0.116) & (0.136) & (0.294) & (0.724) & (3.74) \\
&\multicolumn{9}{l}{Joint} \\
&\quad Grid search & -- & 60 & -- & 0.00160 & 0.01570 & 0.02250 & 0.00170 & 0.00010 \\
& &  &  &  & (0.255) & (0.340) & (0.713) & (2.22) & (21.6) \\
&\quad Profile (Brent) & -- & 100 & -- & 0.00004 & 0.01919 & 0.02642 & 0.00161 & 0.00001 \\
& &  &  &  & (0.0570) & (0.0790) & (0.194) & (1.13) & (8.43) \\
&\quad Nelder--Mead & 100 & 100 & 100 & 0.00004 & 0.01920 & 0.02642 & 0.00162 & 0.00001 \\
& &  &  &  & (0.340) & (0.466) & (1.22) & (5.05) & (36.6) \\
&\quad BFGS & 96 & 66 & 100 & 0.00004 & 0.01919 & 0.02648 & 0.00160 & 0.00001 \\
& &  &  &  & (0.539) & (0.752) & (1.90) & (7.20) & (49.9) \\
&\quad L-BFGS-B & 100 & 100 & 100 & 0.00004 & 0.01919 & 0.02643 & 0.00162 & 0.00001 \\
& &  &  &  & (0.313) & (0.362) & (0.878) & (1.69) & (13.9) \\
&\quad \texttt{nlminb} & 100 & 100 & 100 & 0.00004 & 0.01919 & 0.02643 & 0.00162 & 0.00001 \\
& &  &  &  & (0.255) & (0.268) & (0.684) & (1.44) & (10.8) \\
\bottomrule
\end{tabular}

\end{table}

\subsection{Comparison with competing estimators}

In this subsection, we compare the FPLD estimator with its special cases, namely the estimators of Menon \citep{Menon1963} and Sadani \citep{Sadani2019}, as well as the closed-form estimators based on L-moments (LM) \citep{Hosking1990} and percentiles \citep{Dubey1967}. We also consider as competitors the iterative method-of-moments (MM) and maximum likelihood (ML) estimators \citep{Teimouri2013}, the modified maximum likelihood (MML) estimator \citep{Jokiel2024}, and Firth's bias-corrected estimator \citep{Firth1993} (see also \citealp{Silva2025} and the references therein).

\begin{figure}[H]
    \centering
\includegraphics[width=1.0\linewidth]{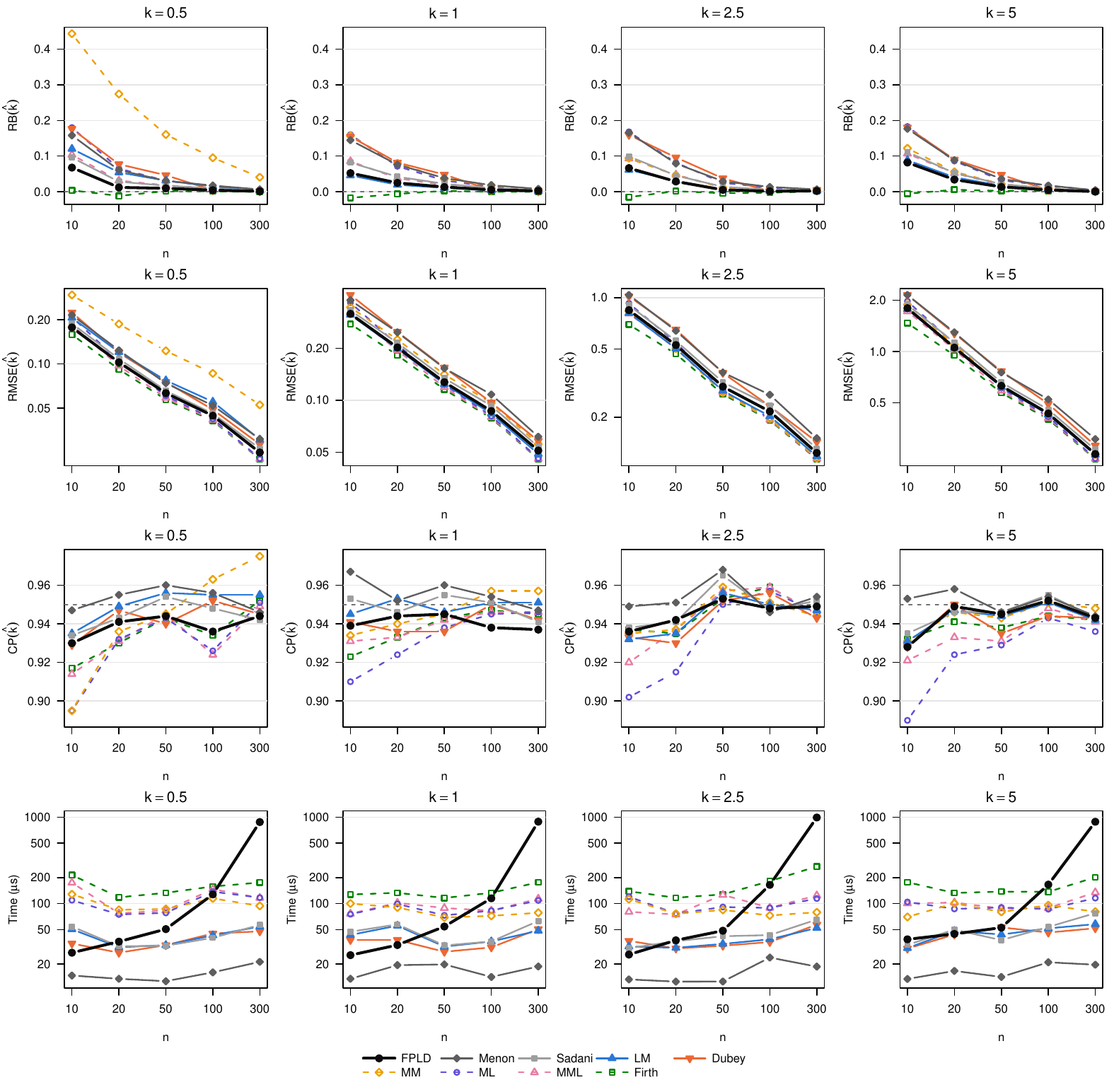}
    \caption{Monte Carlo results for the RB, RMSE, CP, and computational time (in $\mu$s) of the estimators of $k$, with increasing sample size $n \in\{ 10, 20, 50, 100, 300\}$, varying $k \in \{0.5, 1, 2.5, 5\}$, and $\lambda=1$.}
    \label{fig:comp_shape}
\end{figure}

\begin{figure}[H]
    \centering
\includegraphics[width=1.0\linewidth]{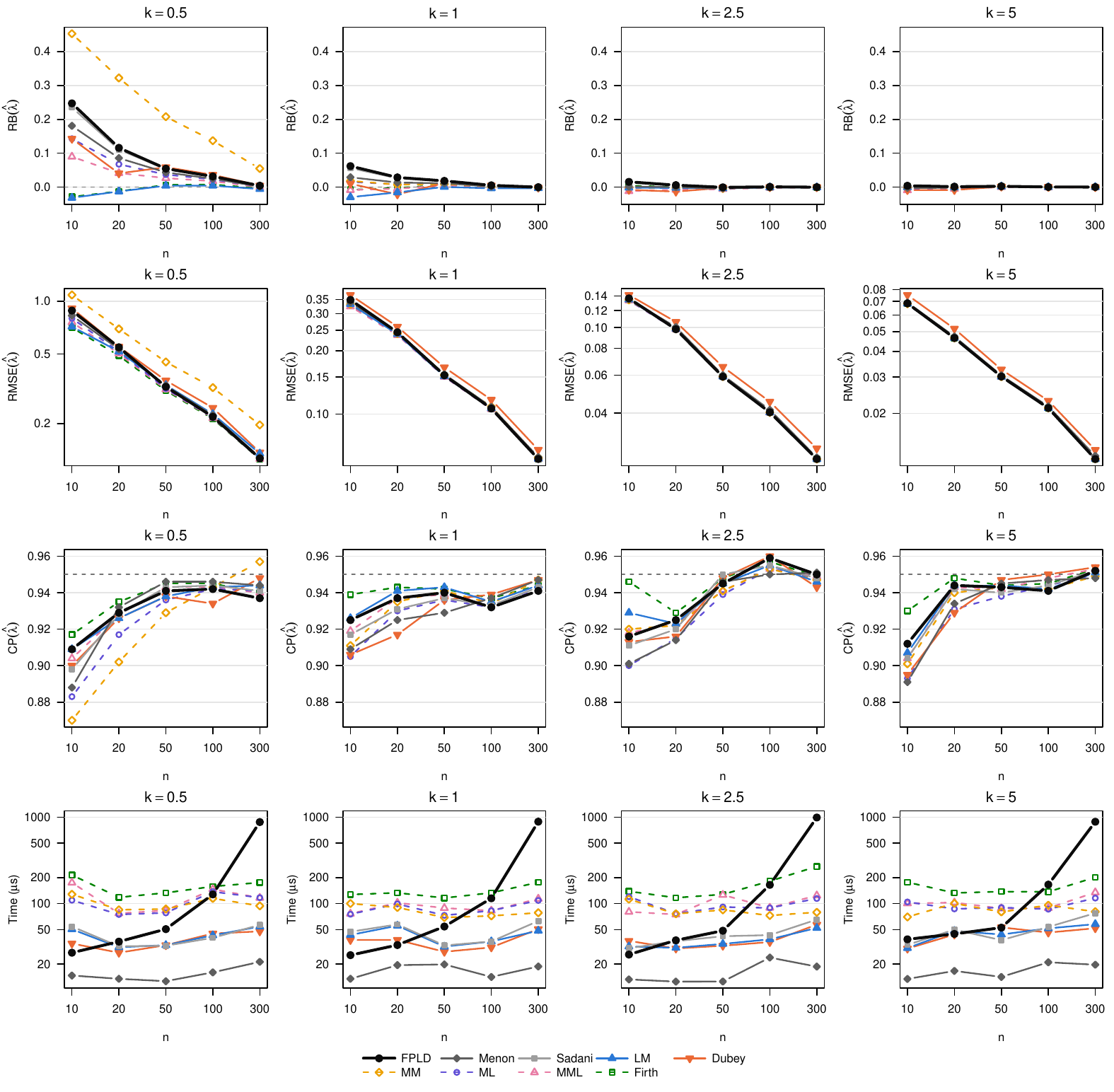}
    \caption{Monte Carlo results for the RB, RMSE, CP, and computational time (in $\mu$s) of the estimators of $\lambda=1$, with increasing sample size $n \in\{ 10, 20, 50, 100, 300\}$ and varying $k \in \{0.5, 1, 2.5, 5\}$.}
    \label{fig:comp_scale}
\end{figure}

\begin{figure}[H]
    \centering
\includegraphics[width=1.0\linewidth]{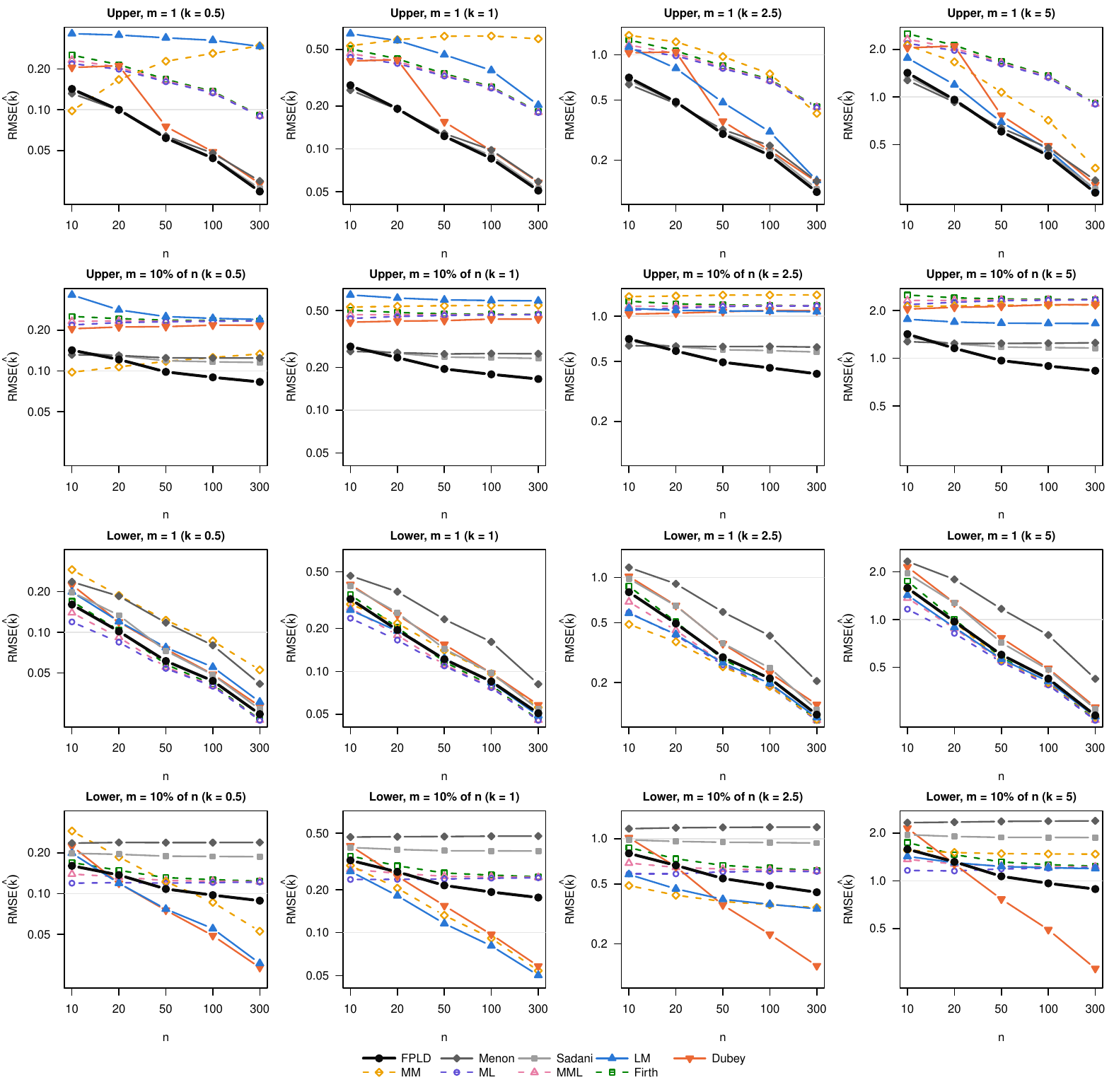}
    \caption{Monte Carlo results for the RMSE of the estimators of $k$ under upper and lower contamination of $m=1$ or 10\% of the observations, with increasing sample size $n \in \{10, 20, 50, 100, 300\}$, varying $k \in \{0.5, 1, 2.5, 5\}$, and $\lambda=1$.}
    \label{fig:robust_shape}
\end{figure}

\begin{figure}[H]
    \centering
\includegraphics[width=1.0\linewidth]{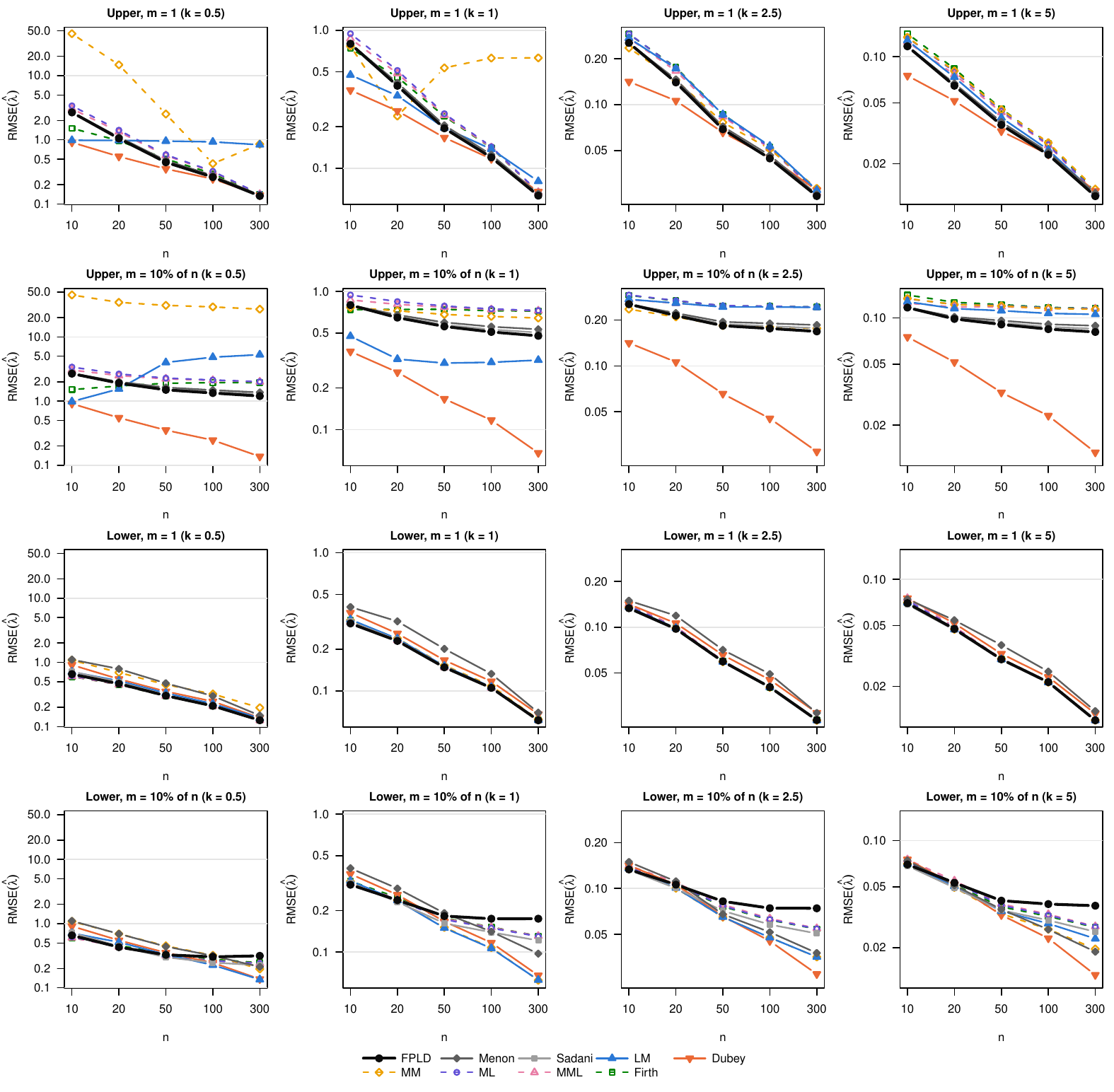}
    \caption{Monte Carlo results for the RMSE of the estimators of $\lambda=1$ under upper and lower contamination of $m=1$ or 10\% of the observations, with increasing sample size $n \in\{ 10, 20, 50, 100, 300\}$ and varying $k \in \{0.5, 1, 2.5, 5\}$.}
    \label{fig:robust_scale}
\end{figure}

Besides RB, RMSE, and CP, we also compare the computational time of the estimators (in $\mu$s). Figures~\ref{fig:comp_shape} and~\ref{fig:comp_scale} show the results for the estimation of $k$ and $\lambda$, respectively. Figures~\ref{fig:robust_shape} and~\ref{fig:robust_scale} show, for $k$ and $\lambda$ respectively, the results in the presence of upper and lower contamination, that is, when the largest values are increased, or the smallest ones are brought close to zero. The RMSE measures the impact on the estimation of the true parameters. In Figures~\ref{fig:comp_shape}--\ref{fig:robust_scale}, FPLD uses $(r,a)$ and $c_r$ from Table~\ref{tab:lookup}.

For $k$ estimation, it is the most accurate closed-form estimator when $k=0.5$, and it outperforms ML for $n\leq 20$, with Wald coverage between 0.93 and 0.95. Under contamination, FPLD is the most robust estimator of $k$ against large values, but not against small ones. 

For $\lambda$ estimation, the estimators are similar except for $k=0.5$ and small $n$, where LM and Firth are better. Requiring no optimization, FPLD is faster than the iterative estimators for $n\leq 20$, but its cost grows as $n$ increases. 
Under contamination, Dubey's is the least affected.
Hence, as a closed-form estimator with no optimization, FPLD rivals ML and beats it for $k$ with small $n$ or large outliers.

Therefore, FPLD is the most accurate closed-form estimator of $k$ (uniformly in $k$), outperforms ML in small samples, yields Wald intervals with near-nominal coverage from $n = 10$, and is robust to atypical observations in both tails, at the cost of a moderate efficiency loss in large samples.

\section{Concluding remarks}
\label{Concluding remarks}

In this paper, we proposed a class of explicit estimators for the shape and scale parameters of the two-parameter Weibull distribution, referred to as the FPLD estimators. Their construction is based on fractional absolute moments of pairwise logarithmic differences, together with weighted pairwise logarithmic minima and maxima.

Strong consistency and joint asymptotic normality have been established. In addition, the determinant-based optimization problem has been reduced analytically to a
one-dimensional problem. An extensive numerical evaluation of the FPLD estimator was conducted, examining the impact of the fractional order, proposing a reference table for selecting the tuning parameters, comparing this approach with optimization-based tuning procedures, and evaluating the performance of the FPLD estimator against competitors.

Several questions remain open. These include the existence and uniqueness of the optimal fractional order, the numerical evaluation of the optimal tuning parameters, and the finite-sample behavior of the resulting estimators. These issues will determine whether the proposed class provides a meaningful practical improvement over
existing explicit estimators.

%\section*{Acknowledgment}
%CNPq, FAPDF, DPI, PPGE/UnB

\vspace*{-0.15cm}

\section*{Declarations}

\paragraph*{Ethics Approval}
Not applicable.

% \paragraph*{Funding Declaration}
% This study was financed in part by CAPES (Finance Code 001). 

%\paragraph*{Acknowledgements}
%The research was supported in part by CNPq and CAPES grants from the Brazilian government.

\paragraph*{Disclosure statement}
There are no conflicts of interest to disclose.

\paragraph*{Author Contributions Statement}
All authors contributed equally to the conception and design of the study, data analysis, interpretation of the results, manuscript preparation, and revision. All authors read and approved the final version of the manuscript.

\paragraph*{Data availability}
% The data used in this study are publicly available from the ... in Data at \url{...}. 
No data were used in this study.

%\vspace*{-0.15cm}
{

}

%\bibliographystyle{apalike}
%\bibliography{refbib}

\appendix
\section{Appendix}\label{ap1}

In this appendix, we provide a reference table for selecting the tuning parameters $(r,a)$ of the FPLD estimators according to the sample size $n$.
We also precomputed 
$$c_r=2\Gamma(r+1)\sum_{j\geqslant 1}\frac{(-1)^{j-1}}{j^r},$$ 
which allows \eqref{eq:kestimator} to be rewritten as
$$\widehat{k}_r=\left(\frac{c_r}{U_{n,r}}\right)^{1/r}.$$

% Lookup table of the tuning parameters of the FPLD estimators. For each sample size $n$: the pair $(r,a)$ chosen by the two-step rule and the constant $c_r=2\Gamma(r+1)\sum_{j\geqslant 1}(-1)^{j-1}/j^r$ in \eqref{eq:logisticmoment}, so that $\widehat{k}_r=(c_r/U_{n,r})^{1/r}$ in \eqref{eq:kestimator}. The values do not depend on $(k,\lambda)$ or on the data. They were computed from $M=20\,000$ standard Weibull samples ($k=\lambda=1$), common to all $n$: $r$ minimizes $n\,\mathrm{MSE}(\log\widehat{k}_r)$ and $a$ minimizes $n\,\mathrm{MSE}(\log\widehat{\lambda}_{a,r})$ for that $r$, both rounded to two decimals. Both criteria are flat near their minima: changing $r$ by $0.02$ or $a$ by $0.05$ increases them by less than 0.1\%. For $n>300$, use $r=0.01$ and $a=0.21$; between two tabulated sample sizes, $r$ and $a$ can be interpolated linearly on $\log n$.

\begin{table}[H]
\centering
\caption{Reference table.}
\label{tab:lookup}
\footnotesize
\setlength{\tabcolsep}{4pt}
\begin{tabular}{rrrr@{\hspace{2.2em}}rrrr@{\hspace{2.2em}}rrrr}
\toprule
$n$ & $r$ & $a$ & $c_r$ & $n$ & $r$ & $a$ & $c_r$ & $n$ & $r$ & $a$ & $c_r$ \\
\cmidrule(lr){1-4}\cmidrule(lr){5-8}\cmidrule(lr){9-12}
5 & 0.51 & 0.26 & 1.0760 & 19 & 0.26 & 0.21 & 1.0067 & 45 & 0.12 & 0.18 & 0.9939 \\
6 & 0.49 & 0.29 & 1.0684 & 20 & 0.26 & 0.21 & 1.0067 & 50 & 0.11 & 0.18 & 0.9937 \\
7 & 0.47 & 0.26 & 1.0611 & 21 & 0.26 & 0.22 & 1.0067 & 60 & 0.10 & 0.21 & 0.9937 \\
8 & 0.44 & 0.26 & 1.0509 & 22 & 0.25 & 0.21 & 1.0052 & 70 & 0.09 & 0.21 & 0.9938 \\
9 & 0.42 & 0.27 & 1.0445 & 23 & 0.24 & 0.21 & 1.0037 & 80 & 0.08 & 0.21 & 0.9940 \\
10 & 0.40 & 0.24 & 1.0385 & 24 & 0.24 & 0.22 & 1.0037 & 90 & 0.06 & 0.20 & 0.9948 \\
11 & 0.38 & 0.27 & 1.0329 & 25 & 0.23 & 0.22 & 1.0023 & 100 & 0.05 & 0.19 & 0.9953 \\
12 & 0.36 & 0.25 & 1.0276 & 26 & 0.22 & 0.22 & 1.0011 & 120 & 0.04 & 0.19 & 0.9960 \\
13 & 0.34 & 0.24 & 1.0227 & 27 & 0.21 & 0.21 & 0.9999 & 140 & 0.01 & 0.19 & 0.9988 \\
14 & 0.32 & 0.24 & 1.0181 & 28 & 0.20 & 0.21 & 0.9988 & 160 & 0.01 & 0.19 & 0.9988 \\
15 & 0.30 & 0.21 & 1.0140 & 29 & 0.20 & 0.21 & 0.9988 & 180 & 0.01 & 0.20 & 0.9988 \\
16 & 0.29 & 0.21 & 1.0120 & 30 & 0.19 & 0.20 & 0.9978 & 200 & 0.01 & 0.21 & 0.9988 \\
17 & 0.28 & 0.21 & 1.0102 & 35 & 0.17 & 0.20 & 0.9962 & 250 & 0.01 & 0.21 & 0.9988 \\
18 & 0.27 & 0.21 & 1.0084 & 40 & 0.15 & 0.19 & 0.9949 & 300 & 0.01 & 0.21 & 0.9988 \\
\bottomrule
\end{tabular}
\end{table}

\end{document}